\documentclass[12pt]{article}     
\usepackage[utf8]{inputenc}

\usepackage[T1]{fontenc}

\usepackage[warn]{mathtext}          
\usepackage[T2A]{fontenc}            

\usepackage{graphicx}      
\usepackage{psfrag}        
\usepackage{caption2}      
\usepackage{soul}          
\usepackage{fancyhdr}      
\usepackage{multirow}      
\usepackage{ltxtable}      
\usepackage{paralist}      
\usepackage[perpage]{footmisc} 

\usepackage{amsmath}
\usepackage{amsfonts}
\usepackage{amssymb}
\usepackage{amsthm}
\usepackage{graphicx,color,framed} 

\usepackage{tikz}

\DeclareMathOperator{\Res}{Re}
\DeclareMathOperator{\Ims}{Im}

\DeclareMathOperator{\mmod}{mod}
\newcommand{\Z}{\scriptstyle}
\newcommand{\D}{\displaystyle}

\renewcommand{\le}{\operatorname{\leqslant}}
\renewcommand{\ge}{\operatorname{\geqslant}}

\usepackage[a4paper, top=30mm, left=25mm, right=25mm, bottom=30mm]{geometry}

\newtheorem{lemma}{Lemma}
\newtheorem{theorem}{Theorem}

\newenvironment{myproof2} {\par\noindent{\bf Proof.}} {\hfill$\scriptstyle\blacksquare$}

\begin{document}

\begin{center}
\textsc{\textbf{\Large On the asymptotic formulae for some multiplicative functions in short intervals}}
\footnote{This research was supported by the grant of Russian Fund of Fundamental Researches  12-01-31165.}
\end{center}

\begin{flushleft}
\begin{center}A.A. Sedunova\\ \end{center}
\end{flushleft}

{\bf Abstract.} In this paper we use the method proposed in \cite{me} in
order to find the mean values of some multiplicative functions connected
with the divisor function on the short interval of summation. Further
investigations have shown that the same technique gives the result for
$$f(n) = \frac{1}{\tau(n^2)}, \;\frac{1}{(\tau(n))^2}, \; \frac{1}{2^{\omega(n)}}, \;\frac{1}{2^{\Omega(n)}}.$$

\begin{section}{Introduction}
	
	In 1919, S.Ramanujan \cite{rmnj} announced the formula
	
	\begin{equation} \label{rmnj}
		\sum_{n \le X}{\frac{1}{\tau(n)}} = {\frac {X}{\sqrt{\ln X}} \left( A_0 + \frac{A_1}{\ln X} + \frac{A_2}{(\ln X)^2} + \ldots + \frac{A_N}{(\ln X)^N} + O \left( \frac{1}{(\ln X)^{N+1}}\right)\right)},
	\end{equation}
	where $A_j$ are some constants,
	$$A_0 = \frac{1}{\sqrt{\pi}} \prod_{p} { \sqrt{p (p-1)}\; \ln{\frac{p}{p-1}}},$$
	$\tau(n)$ denotes the number of divisors of $n$ and $N \ge 0$ is a fixed integer.	
	The complete proof of (\ref{rmnj}) was published in 1922 by B.M.Wilson \cite{wilson}.
	
	In this paper we generalize (\ref{rmnj}) and some other theorems of this type to the case when $n$ runs through the short interval of summation,
	i.e. the interval $x < n \le x+h$, where $x \rightarrow +\infty$ and $h \ll x^{\alpha}$ for a fixed $\alpha$, $0 < \alpha <1$.

	Suppose that $k \ge 2$ is fixed. We use the following definitions
	$$\omega(n) = \omega(p_1^{r_1} p_2^{r_2} \ldots p_s^{r_s}) = s, \;\;
	\Omega(n) =  \Omega(p_1^{r_1} p_2^{r_2} \ldots p_s^{r_s}) = r_1 + r_2 + \ldots + r_s,$$ 
	where $p_1, p_2, \ldots $ are distinct prime numbers, $r_1, r_2, \ldots, r_s$ are positive integers.
	Let us define the multiplicative functions $f_j(n)$, $j=1,2,3,4,$ by the following relations:
	$$f_1(n) = \frac{1}{\tau(n^2)};\quad f_2(n) =  \D {\frac{\D 1}{\D (\tau(n))^2}};\quad f_3(n) =\frac{1}{2^{\omega(n)}}; \quad f_4(n) = \frac{1}{2^{\Omega(n)}}.$$
	Finally, let
	$$S_j(x;h) = \sum_{x < n \le x+h}{f_j(n)}.$$
	Our goal is to prove the following theorems.
	
	\begin{theorem} \label{theorem1}
		The formula
		$$S_1(x;h) = \sum_{x < n \le x+h}{\frac{1}{\tau(n^2)}} = \frac{h}{\left(\ln x\right)^{\frac{2}{3}}} \left( \sum_{n=0}^{n=N}{\frac{A_n}{(\ln x)^n}} + O \left(\frac {1}{(\ln x)^{N+1}} \right)\right),$$
		holds true for any fixed $N \ge 0$ and for $h$ under the following conditions:
		$x^{\alpha_1} e^{(\ln x)^{0,1}} \le h \le x,$ 
		where $\alpha_1=\frac{185}{308}$
		and $A_n$ denote some positive constants that depend only on $n$. 
	\end{theorem}
	
	\begin{theorem} \label{theorem2}
		The formula
		\begin{equation*}
		S_2(x; h) = \sum_{x < n \le x+h}{\frac{1}{(\tau(n))^2}} = \frac{h}{(\ln x)^{\frac{3}{4}}} \left( \sum_{n=0}^{n=N}{\frac{B_n}{(\ln x)^n}} + O \left(\frac {1}{(\ln x)^{N+1}}\right)\right)
		\end{equation*}
		holds true for any fixed $N \ge 0$ and for $h$ under the following conditions:
		$x^{\alpha_2} e^{(\ln x)^{0,1}} \le h \le x,$ 
		where $\alpha_2=\frac{303}{508}$
		and $B_n$ denote some positive constants that depend only on $n$.
	\end{theorem}
	
	\begin{theorem} \label{theorem3}
			The formula
		\begin{equation*}
		S_3(x; h) = \sum_{x < n \le x+h}{\frac{1}{2^{\omega(n)}}} = \frac{h}{(\ln x)^{\frac{1}{2}}} \left( \sum_{n=0}^{n=N}{\frac{C_n}{(\ln x)^n}} + O \left(\frac {1}{(\ln 					x)^{N+1}}\right)\right)
		\end{equation*}
		holds true for any fixed $N \ge 0$ and for $h$ under the following conditions:
		$x^{\alpha_3} e^{(\ln x)^{0,1}} \le h \le x,$ 
		where $\alpha_3=\frac{319}{524}$
		and $C_n$ denote some positive constants that depend only on $n$.
	\end{theorem}
	
	\begin{theorem} \label{theorem4}
			The formula
		\begin{equation*}
		S_4(x; h) = \sum_{x < n \le x+h}{\frac{1}{2^{\Omega(n)}}} = \frac{h}{(\ln x)^{\frac{1}{2}}} \left( \sum_{n=0}^{n=N}{\frac{D_n}{(\ln x)^n}} + O \left(\frac {1}{(\ln 					x)^{N+1}}\right)\right)
		\end{equation*}
		holds true for any fixed $N \ge 0$ and for $h$ under the following conditions:
		$x^{\alpha_4} e^{(\ln x)^{0,1}} \le h \le x,$ 
		where $\alpha_4=\frac{319}{524}$
		and $D_n$ denote some positive constants that depend only on $n$.
	\end{theorem}
	
	{\bf Notations} \\
	In what follows, $C ,C_1, C_2\ldots$ denote positive absolute constants, which are, generally speaking, different in different relations. 
	The symbol $(a,b)$ stands for the greatest common divisor of integer $a$ and $b$.
	Finally, $\theta, \theta_1, \theta_2, \ldots$ denote complex numbers with absolute values not greater than one, which are different in different relations.
\end{section}

\begin{section}{Auxilliary statements}
We need some auxilliary lemmas in order to prove theorems \ref{theorem1} - \ref{theorem4}.

\begin{lemma}
	Let $p$ be a prime number and let $\alpha \ge 1$. Then
	$$	\tau(p^{\alpha}) = \alpha+1,\;\;
	\omega(p^\alpha) = 1,\;\;
	\Omega(p^\alpha)= \alpha.
	$$
\end{lemma}
	
\begin{lemma}\label{perron}{\bf{(Perron's formula).}} 
	Suppose that the series $f(s) = \sum_{n=1}^{\infty}{a_n n^{-s}}$ converges absolutely for $\sigma > 1$, $|a_n| \le A(n)$, where
	$A(n)$ is a positive monotonicially increasing function of $n$ and
	$$\sum_{n=1}^{\infty}{|a_n| n^{-\sigma}} = O\left((\sigma-1)^{-\alpha}\right)$$
	for some $\alpha > 0$, as $\sigma \to  1+0$.
	Then the formula
	$$\sum_{n \le x}{a_n} = \frac{1}{2\pi i} \int_{b-iT}^{b+iT}{f(s) \frac{x^s}{s} ds} + O \left( \frac{x^b}{T(b-1)^{\alpha}}\right) + O \left( \frac{xA(2x)\ln x}{T}\right)$$
	holds true for any $b$, $1<b \le b_0$, $T \ge 2$, $x = N+\frac{1}{2}$ (the constants in O-symbols depend on $b_0$).
\end{lemma}
For the proof of the lemma, see \cite{vorkara}, pp. 334-336.

\begin{lemma}{\label{evzetasqr}}
	The estimate
	$$\int_{0}^{T}{\left|\zeta\left(\tfrac{1}{2} +it\right)\right|^2 dt} \ll T \ln T$$
	holds true for any $T \ge T_0 >1$.
\end{lemma}

\begin{myproof2}
	This lemma follows immediately from the theorem of Hardy and Littlewood (see, for example, \cite{titchmarch}, pp. 140-142).
\end{myproof2}

\begin{lemma} \label{evnofzeta}
	Let $\rho(u) = \frac{1}{2} - \{u\}$. Then the formula
	$$\zeta(s) = \frac{1}{2} + \frac{1}{s-1} +s\int_{1}^{\infty}{\frac{\rho(u) du}{u^{s+1}}},$$
	holds true for $s \neq 1$, $\Res s > 0$.
\end{lemma}
	For the proof, see \cite{vorkara}, pp. 24-25.

\begin{lemma}
For $|t| \ge t_0 >1$ and $\frac{\Z 1}{\Z 2} \le \sigma \le 1+ \frac{1}{\ln t}$ we have
$$|\zeta(\sigma + it)| \ll t^{\frac{c(1 - \sigma)}{2}} \ln t,\;\; \left| L(\sigma + it, \chi_4) \right| \ll t^{\frac{c(1 - \sigma)}{2}} \ln t,$$
where $c = \frac{64}{205}$.
\end{lemma}
It was proved by M.N. Huxley in \cite{huxley2}.

\begin{lemma} \label{T^()lnT}
	The estimates
	$$|\zeta(\sigma +it)| \ll t^{c(1 - \sigma)} \ln t,\;\; \left| L(\sigma +it, \chi_4) \right| \ll t^{c(1 - \sigma)} \ln t$$
	hold true for $|t| \ge t_0 >1$, $c = \frac{64}{205}$ and $\frac{\Z 1}{\Z 2} \le \sigma \le 1+ \frac{1}{\ln t}$.
\end{lemma}
\begin{myproof2}
	Let us consider the entire function defined by the relation
	$$\varphi(s) = \zeta(s) (1-s)$$
	and let $\frac{1}{2} = \sigma_1 \le \sigma \le 1 = \sigma_2$.
	Using lemma \ref{T^()lnT}  we state
	$$\left| \varphi\left( \frac{1}{2} + it\right)\right| \ll t^{1 + \frac{c}{2}} \ln t,$$
	$$\left| \varphi\left( 1 + it\right)\right| \ll t \ln t.$$
	Applying the maximum principle (see for example \cite[V]{fcttheory}), we find
	$$\left| \varphi\left( \sigma + it\right)\right| \ll t^{\Theta(\sigma)} \ln t,$$
	where $$\Theta\left(\frac{1}{2}\right) = 1 + \frac{c}{2}, \quad \Theta (1) =1.$$
	Since $\Theta(\sigma)$ is linear function by \cite[V]{fcttheory}, we conclude that
	$$\left| \varphi\left( \sigma + it\right)\right| \ll t^{c(1-\sigma)+1} \ln t.$$
\end{myproof2}

\vspace{0,2cm}

Let $N(\sigma, T)$ be the number of zeros of $\zeta(s)$ in the region $\Res s \ge \sigma, |\Ims s| \le T$.
Suppose that $q \ge 3$ is an integer and let $\chi$ be the Dirichlet's character modulo $q$.
Then the symbol $N(\sigma, T; \chi)$ stands for the number of zeros of the function $L(s, \chi)$ in the same domain.

\begin{lemma} \label{density}
	The estimates
	$$N(\sigma, T) \ll T^{\frac{12}{5} (1-\sigma)} (\ln T)^{44},$$
	$$\sum_{q \le  Q}\; \sideset{}{^*} {\sum}_{\chi \mmod \;Q} {N(\sigma,T; \chi)} \ll (Q^{2}T)^{\frac{12}{5}(1-\sigma)} (\ln QT)^{22}$$
	hold uniformly for $\frac{1}{\Z 2} \le \sigma \le 1$, $T \ge T_0$ and for $Q \ge 2$
	(the symbol $\sum ^*$ means the summation over all primitive characters $\chi$ modulo $q$).
\end{lemma}
For the proof of the first estimate see \cite{huxley}. The second one can be found in \cite{rmch}.

\begin{lemma} \label{vinlim}
	There exist absolute positive constants $t_0$ and $C$ such that $\zeta(s) \neq 0$, $L(s, \chi_4) \neq 0$ in the region
	$$|t| \ge t_0, \;\; \sigma \ge 1 - \varrho(t),\; \varrho(t) = C (\ln \ln t)^{-\frac{1}{3}} (\ln t)^{-\frac{2}{3}}.$$
\end{lemma}
For the proof of this statement for $\zeta(s)$ see \cite[pp. 116-117]{vorkara}. The proof for $L(s, \chi_4)$ can be done by analogy.
\end{section}

\begin{section}{Proof of the main results}
In this section we give the proofs of theorem 1, 2 and 3, 4.
	\begin{subsection}{The mean value of the function $\frac{\D 1}{\D \tau(n^2)}$ on the short interval}
		Suppose that $\sigma = \Res s > 1$ and let
		$$F(s) = \sum_{n=1}^{\infty}{\frac{1}{\tau(n^2)} \cdot n^{-s}}.$$
		This series converges absolutely, since 
		$$|F(s)| \le \sum_{n=1}^{\infty}{\left|\frac{1}{\tau(n^2)}\right| \cdot n^{-\sigma}} \le \frac{1}{3} \sum_{n=1}^{\infty}{n^{-\sigma}} = 
		\frac{1}{3} \left( 1 + \int_{1}^{\infty}{\frac{du}{u^\sigma}}\right) = \frac{1}{3} \left( 1 + \frac{1}{\sigma-1}\right).$$
		
		Setting $a_n = \displaystyle{\frac{1}{\tau(n^2)}}, A(n) \equiv 1, b = 1 +\displaystyle{\frac{1}{\ln x}}, \alpha =1$ in lemma \ref{perron}, we get
		$$S_1=S(x,h;f_1)=I+O(R),$$
		where
		$$I=\frac{1}{2\pi i}\int_{b-iT}^{b+iT}{F(s) \frac{(x+h)^s-x^s}{s}ds},\;\;\; R=\frac{x^b}{T(b-1)}+\frac{xA(2x)\ln x}{T} \ll \frac{x \ln x}{T}.$$
		Further, $F(s) = \prod_p{F_p(s)}$, where
		$$F_p(s)=1 + \frac{1}{\tau(p^2)p^{s}} + \frac{1}{\tau(p^4)p^{2s}} + \ldots = 1+\frac{1}{3p^s} + \frac{1}{5p^{2s}} + \ldots$$
		Writing $F_p(s)$ in the form
		$$F_p(s)=\left( 1 - \frac{1}{p^s}\right)^{-\frac{1}{3}} \left( 1 - \frac{1}{p^{2s}}\right)^{\frac{1}{45}} G_p(s),$$
		we obtain
		$$F(s)=\frac{(\zeta(s))^{\frac{1}{3}}}{(\zeta(2s))^{\frac{1}{45}}} G(s),$$
		where
		$$G(s)=\prod_p{G_p(s)}=\prod_p{\left( 1 - \frac{1}{p^s}\right)^{-\frac{1}{3}} \left( 1 - \frac{1}{p^{2s}}\right)^{\frac{1}{45}} (1+u(s)+v(s))},$$
		$$u(s)=\frac{1}{3p^s}+\frac{2}{5p^{2s}}, \;v(s)=\frac{1}{9p^{3s}}+\frac{1}{17 p^{4s}}+ \ldots$$
		Now we continue the function $F(s)$ to the left of the line $\Res s =1$.
		Suppose that $\frac{1}{2} \le \sigma \le 1$. Then the following estimates hold true:
		\begin{description}
			\item $$|u(s)| \le \frac{1}{3p^\sigma}\left( 1+\frac{3}{5p^\sigma}\right) \le \left(1 + \frac{3}{5\sqrt{2}}\right)\frac{1}{3p^\sigma} < \frac{1}{2p^\sigma};$$
			\item $$|v(s)| \le \frac{1}{9p^{3\sigma}} \left( 1 + \frac{9}{17p^{\sigma}}+ \ldots \right) \le \frac{1}{9p^{3\sigma}} \left( 1 + \frac{1}{p^\sigma} + \ldots \right) \le \frac{1}{2p^{3\sigma}};$$
			\item $$|u(s)+v(s)| \le \frac{1}{2p^\sigma}+\frac{1}{2p^{3\sigma}} < \frac{3}{4 p^\sigma};$$
			\item $$|u(s) \cdot v(s)| \le \frac{1}{4 p^{4\sigma}} \le \frac{1}{4\sqrt{2}p^{3\sigma}} < \frac{1}{4p^{3\sigma}}.$$
		\end{description}
		
		Now let us consider the expansion
		$$\ln(1+u(s)+v(s))=(u+v)-\frac{1}{2}(u+v)^2+\frac{1}{3}(u+v)^3 - \ldots$$
		Obviously, we have
		$$\left|\frac{1}{3}(u+v)^3 - \frac{1}{4}(u+v)^4+\ldots\right| \le \frac{1}{3}\left(\frac{3}{4 p^\sigma}\right)^3+\frac{1}{4}\left(\frac{3}{4p^\sigma}\right)^4 +\ldots \le $$
		$$\le \frac{1}{3}\left(\frac{3}{4p^\sigma}\right)^3\frac{1}{1-\frac{3}{4p^\sigma}} \le \frac{1}{3} \left( \frac{3}{4}\right)^3 \left( 1 - \frac{3}{4\sqrt{2}}\right)^{-1} \frac{1}{p^{3\sigma}} < \frac{13}{24 p^{3\sigma}}.$$
		Next,
		\begin{equation} \label{uv}
			(u+v)-\frac{1}{2}(u+v)^2=\left(u - \frac{u^2}{2}\right) + \left(v - \frac{v^2}{2}-uv \right).
		\end{equation}
		The second term on the right hand of (\ref{uv}) does not exceed in absolute value
		$$\frac{1}{2 p^{3\sigma}} + \frac{1}{16 \sqrt{2} p^{3\sigma}} < \frac{7}{8p^{3\sigma}}.$$
		Moreover,
		$$u-\frac{u^2}{2} = \frac{1}{3p^s}+\frac{13}{90 p^{2s}} + \frac{2 \theta_1}{p^{3\sigma}}.$$
		Using
		$$ \ln \left( 1 - \frac{1}{p^s}\right) = -\frac{1}{p^s}-\frac{1}{2p^{2s}}- \ldots = -\frac{1}{p^s}-\frac{1}{2p^{2s}} + \frac{7\theta_2}{6}\frac{1}{p^{3\sigma}},$$
		$$\ln \left( 1 - \frac{1}{2p^{2s}}\right) = -\frac{1}{p^{2s}} - \ldots = -\frac{1}{p^{2s}} + \frac{5\theta_3}{4} \frac{1}{p^{3\sigma}},$$
		we finally find
		$$\ln G_p(s) = \ln (1+u(s)+v(s)) +\frac{1}{3} \ln \left( 1 - \frac{1}{p^s}\right) - \frac{1}{45} \ln \left(1-\frac{1}{p^{2s}}\right) = \frac{69 \theta}{p^{3\sigma}}.$$
		
		Finally, for $\frac{1}{2} \le \sigma \le 1$ we get
		$$\left| \sum_{p}{\ln G_p(s)}\right| \le C \sum_{p}{\frac{1}{p^{\frac{3}{2}}}} < C,$$
		$$-C \le \ln |G(s)| \le C, \;\;\; e^{-C} \le |G(s)| \le e^{C}.$$
		
		Let $\Gamma$ be the boundary of the rectangle with the vertices $\frac{\Z 1}{\Z 2} \pm iT, b \pm iT$, where the zeros of $\zeta(s)$ of the form
		$\frac{\Z 1}{\Z 2}+i\gamma$, $|\gamma|<T$, are avoided by the semicircles of the infinitely small radius lying to the right of the line $\Res s = \frac{\Z 1}{\Z 2}$,
		the pole of $\zeta(2s)$ at the point $s = \frac{\Z 1}{\Z 2}$ is avoided by two arcs $\Gamma_1$ and $\Gamma_2$ with the radius $\frac{\Z 1}{\Z \ln x}$,
		and let a horizontal cut be drawn from the critical line inside this rectangle to each zero $\rho = \beta + i\gamma$, $\beta > \frac{\Z 1}{\Z 2}$, $|\gamma| < T$.
		Then the function $F(s)$ is analytic inside $\Gamma$.
		By the Cauchy residue theorem,
		
		$$j_0 = - \sum_{k=1}^{8}{j_k} - \sum_{\rho}{j_{\rho}} = - (j_4+j_5) - \sum_{k \neq 4,5}{j_k} - \sum_{\rho}{j_{\rho}}.
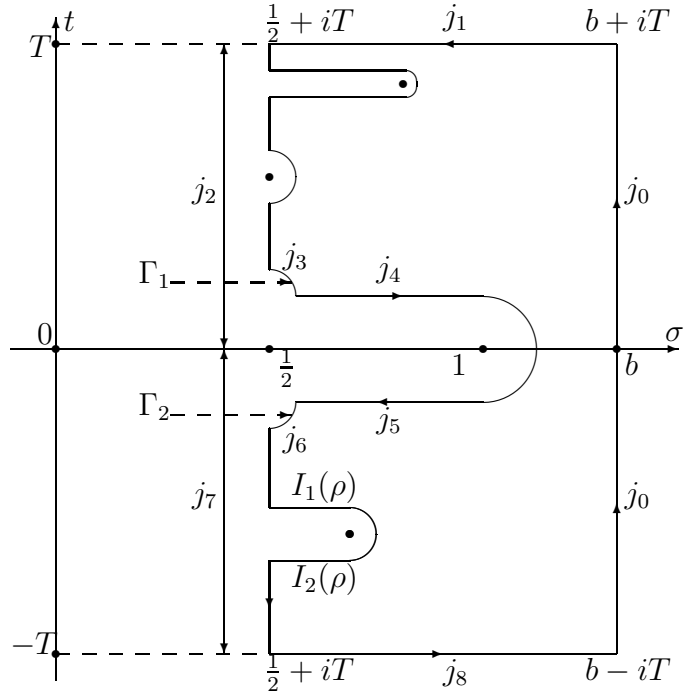
$$
		
		\begin{figure}[tbh]
			\begin{center}
			\begin{picture}(250,267)
			\put(0,125){\vector(1,0){250}}
			\put(245,128){$\sigma$}
			
			\put(17,0){\vector(0,1){250}}
			\put(20,245){$t$}
			
			\put(17,125){\circle*{3}}
			\put(10,127){$0$}
			
			\put(17,10){\circle*{3}}
			\put(0,10){$-T$}
			
			\put(17,240){\circle*{3}}
			\put(7,235){$T$}
			
			\put(97,125){\circle*{3}}
			\put(100,115){$\frac{1}{2}$}
			
			\put(177,125){\circle*{3}}
			\put(165,115){$1$}
			\put(177,125){\oval(40,40)[r]}
			
			\put(97,145){\oval(20,20)[tr]}
			\put(97,105){\oval(20,20)[br]}
			
			\put(127,55){\circle*{3}}
			\put(127,55){\oval(20,20)[r]}
			
			\put(227,125){\circle*{3}}
			\put(230,115){$b$}
			
			\put(227,10){\vector(0,1){57}}
			\put(227,67){\line(0,1){58}}
			\put(230,67){$j_0$}
			
			\put(227,125){\vector(0,1){57}}
			\put(227,182){\line(0,1){58}}
			\put(230,182){$j_0$}
			
			\put(227,240){\vector(-1,0){65}}
			\put(162,240){\line(-1,0){65}}
			\put(162,247){$j_1$}
			
			\put(95,245){$\frac{1}{2}+iT$}
			\put(215,245){$b+iT$}
			
			\put(95,0){$\frac{1}{2}+iT$}
			\put(215,0){$b-iT$}
			
			\put(177,105){\vector(-1,0){40}}
			\put(137,105){\line(-1,0){30}}
			\put(137,95){$j_5$}
			
			\put(97,10){\vector(1,0){65}}
			\put(162,10){\line(1,0){65}}
			\put(162,0){$j_8$}
			
			\put(107,145){\vector(1,0){40}}
			\put(147,145){\line(1,0){30}}
			\put(137,152){$j_4$}
			
			\put(80,182){\vector(0,1){58}}
			\put(80,182){\vector(0,-1){57}}
			\put(68,182){$j_2$}
			
			\put(80,67){\vector(0,1){58}}
			\put(80,67){\vector(0,-1){57}}
			\put(68,67){$j_7$}
			
			\put(102,157){$j_3$}
			\put(102,88){$j_6$}
			
			\put(97,95){\line(0,-1){30}}
			\put(97,45){\vector(0,-1){18}}
			\put(97,27){\line(0,-1){17}}
			
			\put(97,240){\line(0,-1){10}}
			\put(97,220){\line(0,-1){20}}
			
			\put(97,190){\oval(20,20)[r]}
			\put(97,190){\circle*{3}}
			\put(97,180){\line(0,-1){25}}
			
			
			\put(97,65){\line(1,0){30}}
			\put(105,70){$I_1(\rho)$}
			
			\put(97,45){\line(1,0){30}}
			\put(105,35){$I_2(\rho)$}
			
			\put(127,55){\circle*{3}}
			\put(127,55){\oval(20,20)[r]}

			\put(97,230){\line(1,0){50}}
			\put(97,220){\line(1,0){50}}

			\put(147,225){\circle*{3}}
			\put(147,225){\oval(10,10)[r]}
			
			\multiput(17,240)(10,0){8}{\line(1,0){5}}
			\multiput(17,10)(10,0){8}{\line(1,0){5}}
			
			\multiput(60,150)(10,0){5}{\line(1,0){5}}
			\multiput(60,100)(10,0){5}{\line(1,0){5}}
			
			\put(48,150){$\Gamma_1$}
			\put(48,100){$\Gamma_2$}
			
			\put(100,150){\vector(1,0){5}}
			\put(100,100){\vector(1,0){5}}
			
		\end{picture}
			\end{center}
		\end{figure}
		
		By lemma \ref{T^()lnT},
		$$F(s) \ll T^{\frac{1-\sigma}{3}} (\ln T)^{\frac{1}{45}}.$$
		Then
		$$|j_1|=\left| \frac{1}{2 \pi i} \int_{b + iT}^{\frac{1}{2} + iT}{F(s) \frac {(x + h)^s - x^s}{s} ds} \right| \ll
			\frac{1}{T} \int_{\frac{1}{2}}^{b}{T^{\frac{c(1 - \sigma)}{3}} \cdot (\ln x)^{\frac{1}{45}} x^{\sigma} d\sigma} \ll $$
		$$\ll \frac{x}{T} \int _{\frac{1}{2}}^{b}{\frac{x^{\sigma-1}}{T^{\frac{c(\sigma-1)}{3}}} (\ln x)^{\frac{1}{45}}} d\sigma \ll
		\frac{x}{T} \int _{\frac{1}{2}}^{b}{\left(\frac{x}{T^{\frac{c}{3}}} \right)^{\sigma-1} (\ln x)^{\frac{1}{45}}} d\sigma \ll
			\frac{x}{T} (\ln x)^{\frac{1}{45}}.$$
		The similar estimate is valid for $j_8$.
		
		\begin{figure}[tbh]
			\begin{center}
				\begin{picture}(125,150)
			\put(0,75){\vector(1,0){125}}
			\put(121,77){$\sigma$}
			
			\put(10,75){\circle*{3}}
			\put(9,64){$\frac{1}{2}$}
			
			\put(10,85){\oval(40,40)[tr]}
			\put(10,65){\oval(40,40)[br]}
			
			\put(30,85){\vector(1,0){30}}
			\put(60,85){\line(1,0){30}}
			
			\put(10,105){\line(0,1){10}}
			\put(10,35){\line(0,1){10}}
			
			\put(90,65){\vector(-1,0){30}}
			\put(30,65){\line(1,0){30}}
			
			\put(90,75){\oval(20,20)[r]}
			
			\put(25,105){$j_3$}
			\put(25,45){$j_6$}
			
			\put(0,95){$\Gamma_1$}
			\put(0,50){$\Gamma_2$}
			
			\put(12,97){\vector(1,0){13}}
			\put(12,52){\vector(1,0){13}}
			
		\end{picture}
			\end{center}
		\end{figure}
		
		By lemma \ref{evnofzeta}, on $\Gamma_1, \Gamma_2$ we have:
		$$|\zeta(s)| = \left| 0.5 + \frac {1}{s - 1} + s \int_{1}^{\infty}{\frac{\rho(u)}{u^{s + 1}} du} \right| \le 
			0.5 + 2 + 0.1 + 0.6 \cdot 0.5 \int_{1}^{\infty}{\frac {du}{u^{\frac{3}{2}}}} \le 3.2,$$
		$$\left| \zeta(2s)\right| \ge \frac {1}{|2s - 1|} - 1.1 \ge 0.5 \ln x - 1.1 > 0.4 \ln x.$$
		Hence
		$$|F(s)|  \le \left| G(s) \right| \frac {{3.2}^{\frac{1}{3}}}{(0.4 \cdot \ln x)^{\frac{1}{45}}} < C.$$
		Therefore,
		$$|j_3 + j_6| \le \frac {1}{2\pi} \int_{\Gamma_1 \cup \Gamma_2}{|F(s)| \left| \frac {(x + h)^s - x^s}{s} \right| ds} \le
				\frac {C}{2\pi} \int_{- \frac{\pi}{2}}^{\frac{\pi}{2}}{\frac{2 \cdot (2x)^{\frac{1}{2} + \frac{1}{\ln x}}}{\frac{1}{2}} \cdot \frac {d\varphi}{\ln x}} \ll
				\frac{\sqrt{x}}{\ln x}.$$
		Further,
		$$|F(s)| \le |\zeta(s)|^{\frac{1}{3}}(\ln x)^{\frac{1}{45}}|G(s)| \ll (\ln x)^{\frac{1}{45}}\left|\zeta(\sigma+it)\right|^{\frac{1}{3}}.$$
		Hence
		$$|j_2| = \left| \text{p.v.} \frac{1}{2 \pi i}\int_{\frac{1}{2} + iT}^{\frac{1}{2} + \frac{i}{\ln x}}{F(s) \cdot \frac {(x + h)^s - x^s}{s} ds} \right|
			\ll \int_{\frac{1}{\ln x}}^{T} {(\ln x)^{\frac{1}{45}} \cdot \left|\zeta\left(\tfrac{1}{2} + it\right)\right|^{\frac{1}{3}}} \sqrt{x} \frac{dt}{t+1} \ll $$
			$$\ll (\ln x)^{\frac{1}{45}} \sqrt{x} \int_{0}^{T}{\left|\zeta\left(\tfrac{1}{2}+it\right)\right|^{\frac{1}{3}} \frac{dt}{t+1}} =
			(\ln x)^{\frac{1}{45}} \sqrt{x} \sum_{\nu \ge 0}{\int_{T/2^\nu}^{T/2^{\nu +1}}{\frac{|\zeta(\frac{1}{2} +it)|^{\frac{1}{3}}}{t+1}dt}}.$$
		Denoting the summands in the last sum by $j(\nu)$ and taking $X=T\cdot2^{-\nu}$, by the Hölder inequality we get:
		$$j(\nu) \ll \frac{1}{X} \left( \int_{X}^{2X}{\left|\zeta\left(\tfrac{1}{2}+it\right)\right|^2dt}\right)^{\frac{1}{6}} X^{1-\frac{1}{6}} \ll
		\frac{1}{X} \left( X\ln X\right)^{\frac{1}{6}}X^{1-\frac{1}{6}} \ll
		(\ln X)^{\frac{1}{6}} \ll (\ln T)^{\frac{1}{6}}.$$
		Hence,
		$$\sum_{\nu \ge 0}{j(\nu)} \ll (\ln T)^{1+\frac{1}{6}} \ll (\ln T)^{\frac{7}{6}}.$$
		Then the upper bound for $j_2$ has the form
		$$|j_2| \ll (\ln x)^{\frac{1}{45}} \sqrt{x} (\ln x)^{\frac{7}{6}} = \sqrt{x} (\ln x)^{\frac{1}{45}+\frac{7}{6}} \ll \sqrt{x} (\ln x)^{\frac{107}{90}}.$$
		The integral $j_7$ is estimated as above.
		
		The main term arises from the calculation of $j_4$ and $j_5$. Let us define the entire function $w(s)$ by the relation
		$$\zeta(s) = \frac{w(s)}{s-1}$$
		and let $s = 1 - u + i\cdot0$, where $0 \le u \le \frac{1}{2}$. Then
		$$\sqrt[3]{\zeta(s)} = \frac {\sqrt[3]{w(s)}}{\sqrt[3]{-u+i \cdot 0}}.$$
		Since $-u+i\varepsilon \to u\cdot e^{\pi i}$ as $\varepsilon \to +0$, then
		$$\sqrt[3]{-u+i \cdot 0} = \sqrt[3]{u} e^{\frac{\pi i}{3}}, \;\; \sqrt[3]{\zeta(s)}=\frac{\sqrt[3]{w(\sigma)}}{\sqrt[3]{u}} e^{-\frac{\pi i}{3}}.$$
		Therefore, on the upper edge of the cut we have
		$$F(s) = \frac{\sqrt[3]{w(1-u)}}{(\zeta(2-2u))^{\frac{1}{45}}}\; G(1-u)\frac{e^{-\frac{\pi i}{3}}}{\sqrt[3]{u}} = \frac{\Pi(u)e^{-\frac{\pi i}{3}}}{\sqrt[3]{u}},$$
		where
		$$\Pi(u) = G(1-u) \frac{\sqrt[3]{w(1-u)}}{(\zeta(2-2u))^{\frac{1}{45}}}.$$
		
		Hence,
		$$j_4 = \frac{1}{2 \pi i} \int_{\frac{1}{2} + \frac{1}{\ln x} + i \cdot 0}^{1+i \cdot 0}{F(\sigma + i \cdot 0) \frac{(x+h)^s- x^s}{s} ds} = $$
		$$=\frac{1}{2\pi i} \int_{\frac{1}{2}+\frac{1}{\ln x}+i \cdot 0}^{1+i \cdot 0}{F(\sigma+i \cdot 0) \int_{0}^{h}{(x+u)^{s-1}du}\;ds} =$$
		$$= \frac{1}{2\pi i}\int_{x}^{x+h}{\int_{\frac{1}{2}+\frac{1}{\ln x}}^{1}{F(\sigma+i \cdot 0)y^{\sigma-1}d\sigma}\;dy}=
		\frac{e^{-\frac{\pi i}{3}}}{2\pi i} \int_{x}^{x+h}{{\int_{0}^{\frac{1}{2}-\frac{1}{\ln x}}{\frac{\Pi(u)y^{-u}}{\sqrt[3]{u}}}du}\;dy}.$$
		Suppose that $N \ge 0$ is fixed. Then
		$$\Pi(u)= \Pi_0+\Pi_1u + \Pi_2u^2 + \ldots + \Pi_N u^N + O_N(u^{N+1}),$$
		where
		$$\Pi_0=\Pi(0)=\frac{\sqrt[3]{w(1)}}{(\zeta(2))^{\frac{1}{45}}}.$$
		Thus, we have
		$$j_4 = \frac{e^{-\frac{\pi i}{3}}}{2\pi i} \int_{x}^{x+h}{\left(\sum_{0 \le n \le N}{\Pi_n \int_{0}^{\frac{1}{2}-\frac{1}{\ln x}}{\frac{u^{n}y^{-u}}{\sqrt[3]{u}}du}}+ O\left( J\right)\right)dy},$$
		where
		$$J = \int_{0}^{\frac{1}{2}-\frac{1}{\ln x}}{\frac{u^{N+1}y^{-u}}{\sqrt[3]{u}}du} \le  \frac{\Gamma\left(N+2-\frac{1}{3}\right)}{(\ln y)^{N+\frac{5}{3}}}.$$
		Using the estimate
		$$\int_{\lambda}^{\infty}{w^{k-\gamma}e^{-w}dw < ek!\;\lambda^{k-\gamma}e^{-\lambda}},$$
		where $\lambda>1$, $0 < \gamma <1$, $k \ge 1$, we easily get
		$$\int_{0}^{\frac{1}{2}-\frac{1}{\ln x}}{\frac{u^n y^{-u}}{\sqrt[3]{u}}du} =\frac{1}{(\ln y)^{n+\frac{2}{3}}} \int_{0}^{\ln y \left(\frac{1}{2}-\frac{1}{\ln x} \right)}{e^{-w} w^{n-\frac{1}{3}} dw}=$$
		$$=\frac{1}{(\ln y)^{n+\frac{2}{3}}} \left( \int_{0}^{\infty}{e^{-w} w^{n-\frac{1}{3}} dw} + \frac{\theta e n! (\ln y)^{n-\frac{1}{3}}}{\sqrt{y}}\right)=
		\frac{\Gamma\left( n+\frac{2}{3}\right)}{(\ln y)^{n+\frac{2}{3}}}+ \frac{\theta e n!}{\sqrt{y} \ln y}.$$
		Therefore,
		$$j_4 = \frac{e^{-\frac{\pi i}{3}}}{2\pi i} \int_{0}^{h}{\sum_{0 \le n \le N}{\frac{\Pi_n \Gamma(n+\frac{2}{3})}{(\ln (x+u))^{n+\frac{2}{3}}}}du}+O\left( \frac{h}{\left(\ln x\right)^{N+\frac{5}{3}}}\right).$$
		Let
		$$\varphi(x) = \frac{1}{\left( \ln x\right)^{n+\frac{2}{3}}}.$$
		Then the Lagrange mean-value theorem yields
		$$\varphi(x+u) = \varphi(x) + u \varphi^{'}(x+\theta_1 u) = \frac{1}{\left( \ln x\right)^{n+\frac{2}{3}}} + \frac{\theta_2 h \left( n+\frac{2}{3}\right)}{x \left( \ln x\right)^{n+\frac{5}{3}}}.$$
		Thus we get
		$$j_4 = \frac{h e^{-\frac{\pi i}{3}}}{2\pi i} \sum_{0 \le n \le N}{\frac{\Pi_n \Gamma\left( n+\frac{2}{3}\right)}{(\ln x)^{n+\frac{2}{3}}}}+ O\left( \frac{h}{(\ln x)^{N+\frac{5}{3}}}\right) + O\left( \frac{h^2}{x} \frac{1}{(\ln x)^{N+\frac{5}{3}}}\right),$$
		$$j_5 = -\frac{h e^{\frac{\pi i}{3}}}{2\pi i} \sum_{0 \le n \le N}{\frac{\Pi_n \Gamma\left( n+\frac{2}{3}\right)}{(\ln x)^{n+\frac{2}{3}}}}+ O\left( \frac{h}{(\ln x)^{N+\frac{5}{3}}}\right) + O\left( \frac{h^2}{x} \frac{1}{(\ln x)^{N+\frac{5}{3}}}\right).$$
		
		Finally,
		$$-(j_4+j_5)=
		-\frac{h}{(\ln x)^{\frac{2}{3}}} \left( \sum_{0 \le n \le N}{\frac{(-1)^n \Pi_n}{\Gamma \left( \frac{1}{3}-n\right) (\ln x)^n}} + O\left( \frac{1}{(\ln x)^{N+1}}\right) + O \left( \frac{h}{x(\ln x)^{N+1}}\right)\right).$$
		
		It remains to estimate the sum
		$$\sum_{|\gamma|<T}{j_{\rho}},\;j_{\rho}=I_1(\rho)+I_2(\rho),\;\text{where}\;\;\rho = \beta +i\gamma.$$
		
		\begin{figure}[tbh]
			\begin{center}
				\begin{tikzpicture}

					\draw [->](0,3.5) -- (1.25,3.5);
					\draw (1.25,3.5) -- (2.5,3.5);
					\draw (0,2) -- (1.25,2);
					\draw [<-](1.25,2) -- (2.5,2);

					\draw [->](0,2.75) -- (5,2.75);

					\draw (0,4.1) -- (0,3.5);
					\draw (0,1) -- (0,2);

					\fill[black] (3.4,2.75) circle(0.05);

					\draw[black] (2.5,3.5) arc (140:0:1.18);
					\draw[black] (2.5,2) arc (220:360:1.18);

					\draw (1.25,3.5) node [above]  {$I_1(\rho)$};
					\draw (1.25,2) node [below]  {$I_2(\rho)$};
					\draw (3.4,2.75) node [below]  {$\rho=\beta+i\gamma$};
					\draw (5,2.75) node [above]  {$\sigma$};
					\draw (0,1.2) node [right]  {$\scriptstyle{\Res s = \frac{1}{2}}$};

			\end{tikzpicture}
		\end{center}
		\end{figure}
		Since
		$$\left| \frac{(x+h)^s-x^s}{s}\right| = \left| \int_{0}^{h}{(x+u)^{s-1}du}\right| \ll  \int_{0}^{h}{(x+u)^{\sigma-1}du} \ll h x^{\sigma-1},$$
		then
		$$I_1(\rho) \ll \int_{\frac{1}{2}}^{\beta}{(\ln x)^{\frac{1}{45}} |\zeta(\sigma+i\gamma)|^{\frac{1}{3}} h x^{\sigma-1} d\sigma} \ll
		\frac{h}{x} (\ln x)^{\frac{1}{45}} \int_{\frac{1}{2}}^{\beta}{x^{\sigma}  |\zeta(\sigma+i\gamma)|^{\frac{1}{3}} d\sigma}$$
		and the same estimate is valid for $I_2(\rho)$.
		Hence,
		$$|j_\rho| \ll \int_{\frac{1}{2}}^{\beta}{h x^{\sigma-1} (\ln x)^{\frac{1}{45}} T^{\frac{c(1-\sigma)}{3}} (\ln x)^{\frac{1}{3}} d\sigma} \ll
		h (\ln x)^{\frac{1}{45}+\frac{1}{3}} \int_{\frac{1}{2}}^{\beta}{\left( \frac{T^{\frac{c}{3}}}{x}\right)^{1-\sigma} d\sigma} \ll$$
		$$\ll h (\ln x)^{\frac{1}{45}+\frac{1}{3}} \int_{\frac{1}{2}}^{1}{g(\rho,\sigma) \left( \frac{T^{\frac{c}{3}}}{x}\right)^{1-\sigma} d\sigma},$$
		where
		$$g(\rho,\sigma) = \left\{
			\begin{aligned}
				1,\;\text{if}\;\sigma \le \beta,\\ 
				0,\;\text{if}\;\sigma > \beta.\\
			\end{aligned}
			\right.$$
		Applying lemma \ref{density}, we get
		$$\sum_{|\gamma|<T}{j_\rho} \ll h (\ln x)^{\frac{16}{45}} \int_{\frac{1}{2}}^{1-\varrho (T)}{\left( \sum_{|\gamma| < T}{g(\rho; \gamma)} \right) \left( \frac{T^{\frac{c}{3}}}{x}\right)^{1-\sigma} d\sigma} \ll$$
		$$\ll h (\ln x)^{\frac{16}{45}} \int_{\frac{1}{2}}^{1-\varrho (T)} {N(\sigma; T) \left( \frac{T^{\frac{c}{3}}}{x}\right)^{1-\sigma} d\sigma}\ll$$
		$$\ll h (\ln x)^{44+\frac{16}{45}} \int_{\frac{1}{2}}^{1-\varrho (T)} {\left( \frac{T^{\frac{c}{3}}}{x}\right)^{1-\sigma} T^{\frac{12}{5}(1-\sigma)} d\sigma} \ll$$
		$$\ll  h (\ln x)^{45} \int_{\frac{1}{2}}^{1-\varrho (T)} {\left( \frac{T^{\frac{12}{5}+\frac{c}{3}}}{x}\right)^{1-\sigma} d\sigma} \ll
		 h (\ln x)^{45} \left( \frac{T^{\frac{12}{5} + \frac{c}{3}}}{x}\right)^{1-\varrho (T)}.$$
		Let $D(x)=e^{C_1(\ln x)^{0.8}}.$
		Choosing $T$ from the equation,
		$$T^{\frac{12}{5} + \frac{c}{3}}=x D^{-1}(x)$$
		we easily conclude that
		$$T=x^{\left(\frac{c}{3} + \frac{12}{5}\right)^{-1}} D(x)^{-\left(\frac{c}{3} + \frac{12}{5}\right) ^{-1}}.$$
		
		Obviously, the formula for $S_1$ is asymptotic if
		$$h \gg \frac{x}{T} (\ln x)^2 = x^{\alpha_1} e^{C_2 (\ln x)^{0.8}},$$	
		where
		$$\alpha_1 = 1-\frac{123}{308}=\frac{185}{308}.$$
	\end{subsection}

	\begin{subsection}{The mean value of the function $\frac{\D 1}{\D (\tau(n))^2}$ on the short interval}
		Suppose that $\sigma = \Res s > 1$ and let
		$$F(s) = \sum_{n=1}^{\infty}{\frac{1}{(\tau(n))^2} \cdot n^{-s}}.$$
		This series converges absolutely, since 
		$$|F(s)| \le \sum_{n=1}^{\infty}{\left|\frac{1}{(\tau(n))^2}\right| \cdot n^{-\sigma}} \le \frac{1}{4} \sum_{n=1}^{\infty}{n^{-\sigma}} = 
		\frac{1}{4} \left( 1 + \int_{1}^{\infty}{\frac{du}{u^\sigma}}\right) = \frac{1}{4} \left( 1 + \frac{1}{\sigma-1}\right).$$
		
		Setting $a_n = \displaystyle{\frac{1}{(\tau(n))^2}}, A(n) \equiv 1, b = 1 +\displaystyle{\frac{1}{\ln x}}, \alpha =1$ in lemma \ref{perron}, we get
		$$S_2=S(x,h;f_2)=I+O(R),$$
		where
		$$I=\frac{1}{2\pi i}\int_{b-iT}^{b+iT}{F(s) \frac{(x+h)^s-x^s}{s}ds},\;\;\; R=\frac{x^b}{T(b-1)}+\frac{xA(2x)\ln x}{T} \ll \frac{x \ln x}{T}.$$
		Further, $F(s) = \prod_p{F_p(s)}$, where
		$$F_p(s)=1 + \frac{1}{(\tau(p))^2 p^{s}} + \frac{1}{(\tau(p^2))^2p^{2s}} + \ldots = 1+\frac{1}{4 p^s} + \frac{1}{9 p^{2s}} + \ldots$$
		Writing $F_p(s)$ in the form
		$$F_p(s)=\left( 1 - \frac{1}{p^s}\right)^{-\frac{1}{4}} \left( 1 - \frac{1}{p^{2s}}\right)^{\frac{19}{244}} G_p(s),$$
		we obtain
		$$F(s)=\frac{(\zeta(s))^{\frac{1}{4}}}{(\zeta(2s))^{\frac{19}{244}}} G(s),$$
		where
		$$G(s)=\prod_p{G_p(s)}=\prod_p{\left( 1 - \frac{1}{p^s}\right)^{-\frac{1}{4}} \left( 1 - \frac{1}{p^{2s}}\right)^{\frac{19}{244}} (1+u(s)+v(s))},$$
		$$u(s)=\frac{1}{4p^s}+\frac{1}{9p^{2s}}, \;v(s)=\frac{1}{16p^{3s}}+\frac{1}{25 p^{4s}}+ \ldots$$
		Now we continue the function $F(s)$ to the left of the line $\Res s =1$.
		Suppose that $\frac{1}{2} \le \sigma \le 1$. Then the following estimates hold true:
		\begin{description}
			\item $$|u(s)| \le \frac{1}{4p^\sigma}\left( 1+\frac{4}{9p^\sigma}\right) \le \left(1 + \frac{4}{9\sqrt{2}}\right)\frac{1}{4p^\sigma} < \frac{1}{3p^\sigma};$$
			\item $$|v(s)| \le \frac{1}{16p^{3\sigma}} \left( 1 + \frac{16}{25p^{\sigma}}+ \ldots \right) \le \frac{1}{16 p^{3\sigma}} \left( 1 + \frac{1}{p^\sigma} + \ldots \right) \le \frac{1}{4 p^{3\sigma}};$$
			\item $$|u(s)+v(s)| \le \frac{1}{3p^\sigma}+\frac{1}{4p^{3\sigma}} < \frac{1}{2 p^\sigma};$$
			\item $$|u(s) \cdot v(s)| \le \frac{1}{12 p^{4\sigma}} \le \frac{1}{12\sqrt{2}p^{3\sigma}} < \frac{1}{12p^{3\sigma}}.$$
		\end{description}
		
		Now let us consider the expansion
		$$\ln(1+u(s)+v(s))=(u+v)-\frac{1}{2}(u+v)^2+\frac{1}{3}(u+v)^3 - \ldots$$
		Obviously, we have
		$$\left|\frac{1}{3}(u+v)^3 - \frac{1}{4}(u+v)^4+\ldots\right| \le \frac{1}{3}\left(\frac{1}{2 p^\sigma}\right)^3+\frac{1}{4}\left(\frac{1}{2p^\sigma}\right)^4 +\ldots \le $$
		$$\le \frac{1}{3}\left(\frac{1}{2p^\sigma}\right)^3\frac{1}{1-\frac{1}{2p^\sigma}} \le \frac{1}{2} \left( \frac{1}{2}\right)^3 \left( 1 - \frac{1}{2\sqrt{2}}\right)^{-1} \frac{1}{p^{3\sigma}} < \frac{1}{12 p^{3\sigma}}.$$
		Next,
		\begin{equation} \label{uv}
			(u+v)-\frac{1}{2}(u+v)^2=\left(u - \frac{u^2}{2}\right) + \left(v - \frac{v^2}{2}-uv \right).
		\end{equation}
		For the second term on the right hand side we have the upper bound
		$$\frac{1}{4 p^{3\sigma}} + \frac{1}{32 \sqrt{2} p^{3\sigma}} + \frac{1}{12 p^{3 \sigma}} < \frac{3}{8p^{3\sigma}}.$$
		Moreover,
		$$u-\frac{u^2}{2} = \frac{1}{4p^s}+\frac{23}{288 p^{2s}} + \frac{\theta_1}{18 p^{3\sigma}}.$$
		Using
		$$ \ln \left( 1 - \frac{1}{p^s}\right) = -\frac{1}{p^s}-\frac{1}{2p^{2s}}- \ldots = -\frac{1}{p^s}-\frac{1}{2p^{2s}} + \frac{7\theta_2}{6}\frac{1}{p^{3\sigma}},$$
		$$\ln \left( 1 - \frac{1}{2p^{2s}}\right) = -\frac{1}{p^{2s}} - \ldots = -\frac{1}{p^{2s}} + \frac{5\theta_3}{4} \frac{1}{p^{3\sigma}},$$
		we finally find
		$$\ln G_p(s) = \ln (1+u(s)+v(s)) +\frac{1}{4} \ln \left( 1 - \frac{1}{p^s}\right) - \frac{19}{244} \ln \left(1-\frac{1}{p^{2s}}\right) = \frac{15 \theta}{p^{3\sigma}}.$$
		
		Finally, for $\frac{1}{2} \le \sigma \le 1$ we get
		$$\left| \sum_{p}{\ln G_p(s)}\right| \le C \sum_{p}{\frac{1}{p^{\frac{3}{2}}}} < C,$$
		$$-C \le \ln |G(s)| \le C, \;\;\; e^{-C} \le |G(s)| \le e^{C}.$$
		
		Let $\Gamma$ be the boundary of the rectangle with the vertices $\frac{\Z 1}{\Z 2} \pm iT, b \pm iT$, where the zeros of $\zeta(s)$ of the form
		$\frac{\Z 1}{\Z 2}+i\gamma$, $|\gamma|<T$, are avoided by the semicircles of the infinitely small radius lying to the right of the line $\Res s = \frac{\Z 1}{\Z 2}$,
		the pole of $\zeta(2s)$ at the point $s = \frac{\Z 1}{\Z 2}$ is avoided by two arcs $\Gamma_1$ and $\Gamma_2$ with the radius $\frac{\Z 1}{\Z \ln x}$,
		and let a horizontal cut be drawn from the critical line inside this rectangle to each zero $\rho = \beta + i\gamma$, $\beta > \frac{\Z 1}{\Z 2}$, $|\gamma| < T$.
		Then the function $F(s)$ is analytic inside $\Gamma$.
		By the Cauchy residue theorem,
		
		$$j_0 = - \sum_{k=1}^{8}{j_k} - \sum_{\rho}{j_{\rho}} = - (j_4+j_5) - \sum_{k \neq 4,5}{j_k} - \sum_{\rho}{j_{\rho}}.
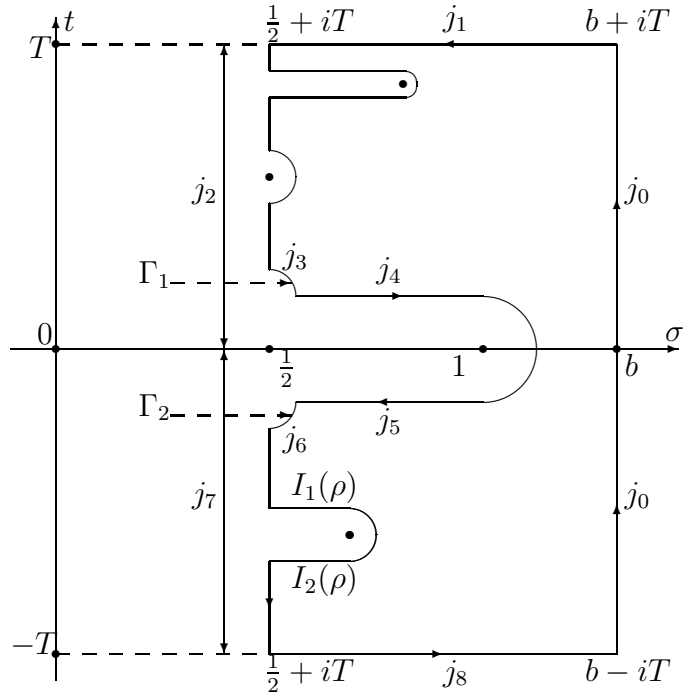
$$
		
		\begin{figure}[tbh]
			\begin{center}
			\begin{picture}(250,267)
			\put(0,125){\vector(1,0){250}}
			\put(245,128){$\sigma$}
			
			\put(17,0){\vector(0,1){250}}
			\put(20,245){$t$}
			
			\put(17,125){\circle*{3}}
			\put(10,127){$0$}
			
			\put(17,10){\circle*{3}}
			\put(0,10){$-T$}
			
			\put(17,240){\circle*{3}}
			\put(7,235){$T$}
			
			\put(97,125){\circle*{3}}
			\put(100,115){$\frac{1}{2}$}
			
			\put(177,125){\circle*{3}}
			\put(165,115){$1$}
			\put(177,125){\oval(40,40)[r]}
			
			\put(97,145){\oval(20,20)[tr]}
			\put(97,105){\oval(20,20)[br]}
			
			\put(127,55){\circle*{3}}
			\put(127,55){\oval(20,20)[r]}
			
			\put(227,125){\circle*{3}}
			\put(230,115){$b$}
			
			\put(227,10){\vector(0,1){57}}
			\put(227,67){\line(0,1){58}}
			\put(230,67){$j_0$}
			
			\put(227,125){\vector(0,1){57}}
			\put(227,182){\line(0,1){58}}
			\put(230,182){$j_0$}
			
			\put(227,240){\vector(-1,0){65}}
			\put(162,240){\line(-1,0){65}}
			\put(162,247){$j_1$}
			
			\put(95,245){$\frac{1}{2}+iT$}
			\put(215,245){$b+iT$}
			
			\put(95,0){$\frac{1}{2}+iT$}
			\put(215,0){$b-iT$}
			
			\put(177,105){\vector(-1,0){40}}
			\put(137,105){\line(-1,0){30}}
			\put(137,95){$j_5$}
			
			\put(97,10){\vector(1,0){65}}
			\put(162,10){\line(1,0){65}}
			\put(162,0){$j_8$}
			
			\put(107,145){\vector(1,0){40}}
			\put(147,145){\line(1,0){30}}
			\put(137,152){$j_4$}
			
			\put(80,182){\vector(0,1){58}}
			\put(80,182){\vector(0,-1){57}}
			\put(68,182){$j_2$}
			
			\put(80,67){\vector(0,1){58}}
			\put(80,67){\vector(0,-1){57}}
			\put(68,67){$j_7$}
			
			\put(102,157){$j_3$}
			\put(102,88){$j_6$}
			
			\put(97,95){\line(0,-1){30}}
			\put(97,45){\vector(0,-1){18}}
			\put(97,27){\line(0,-1){17}}
			
			\put(97,240){\line(0,-1){10}}
			\put(97,220){\line(0,-1){20}}
			
			\put(97,190){\oval(20,20)[r]}
			\put(97,190){\circle*{3}}
			\put(97,180){\line(0,-1){25}}
			
			
			\put(97,65){\line(1,0){30}}
			\put(105,70){$I_1(\rho)$}
			
			\put(97,45){\line(1,0){30}}
			\put(105,35){$I_2(\rho)$}
			
			\put(127,55){\circle*{3}}
			\put(127,55){\oval(20,20)[r]}

			\put(97,230){\line(1,0){50}}
			\put(97,220){\line(1,0){50}}

			\put(147,225){\circle*{3}}
			\put(147,225){\oval(10,10)[r]}
			
			\multiput(17,240)(10,0){8}{\line(1,0){5}}
			\multiput(17,10)(10,0){8}{\line(1,0){5}}
			
			\multiput(60,150)(10,0){5}{\line(1,0){5}}
			\multiput(60,100)(10,0){5}{\line(1,0){5}}
			
			\put(48,150){$\Gamma_1$}
			\put(48,100){$\Gamma_2$}
			
			\put(100,150){\vector(1,0){5}}
			\put(100,100){\vector(1,0){5}}
			
		\end{picture}
			\end{center}
		\end{figure}
		
		By lemma \ref{T^()lnT},
		$$F(s) \ll T^{\frac{c(1-\sigma)}{4}} (\ln T)^{\frac{19}{244}}.$$
		Then
		$$|j_1|=\left| \frac{1}{2 \pi i} \int_{b + iT}^{\frac{1}{2} + iT}{F(s) \frac {(x + h)^s - x^s}{s} ds} \right| \ll
			\frac{1}{T} \int_{\frac{1}{2}}^{b}{T^{\frac{c(1 - \sigma)}{4}} \cdot (\ln x)^{\frac{19}{244}} x^{\sigma} d\sigma} \ll $$
		$$\ll \frac{x}{T} \int _{\frac{1}{2}}^{b}{\frac{x^{\sigma-1}}{T^{\frac{c(\sigma-1)}{4}}} (\ln x)^{\frac{19}{244}}} d\sigma \ll
		\frac{x}{T} \int _{\frac{1}{2}}^{b}{\left(\frac{x}{T^{\frac{c}{4}}} \right)^{\sigma-1} (\ln x)^{\frac{19}{244}}} d\sigma \ll
			\frac{x}{T} (\ln x)^{\frac{19}{244}}.$$
		The similar estimate is valid for $j_8$.
		
		\begin{figure}[tbh]
			\begin{center}
				\begin{picture}(125,150)
			\put(0,75){\vector(1,0){125}}
			\put(121,77){$\sigma$}
			
			\put(10,75){\circle*{3}}
			\put(9,64){$\frac{1}{2}$}
			
			\put(10,85){\oval(40,40)[tr]}
			\put(10,65){\oval(40,40)[br]}
			
			\put(30,85){\vector(1,0){30}}
			\put(60,85){\line(1,0){30}}
			
			\put(10,105){\line(0,1){10}}
			\put(10,35){\line(0,1){10}}
			
			\put(90,65){\vector(-1,0){30}}
			\put(30,65){\line(1,0){30}}
			
			\put(90,75){\oval(20,20)[r]}
			
			\put(25,105){$j_3$}
			\put(25,45){$j_6$}
			
			\put(0,95){$\Gamma_1$}
			\put(0,50){$\Gamma_2$}
			
			\put(12,97){\vector(1,0){13}}
			\put(12,52){\vector(1,0){13}}
			
		\end{picture}
			\end{center}
		\end{figure}
		
		Using the estimations from theorem \ref{theorem1}, on $\Gamma_1, \Gamma_2$ we have:
		$$|\zeta(s)| \le 3.2, \quad \left| \zeta(2s)\right| > 0.4 \ln x.$$
		Hence
		$$|F(s)|  \le \left| G(s) \right| \frac {{3.2}^{\frac{1}{4}}}{(0.4 \cdot \ln x)^{\frac{19}{244}}} < C.$$
		Therefore,
		$$|j_3 + j_6| \le \frac {1}{2\pi} \int_{\Gamma_1 \cup \Gamma_2}{|F(s)| \left| \frac {(x + h)^s - x^s}{s} \right| ds} \le
				\frac {C}{2\pi} \int_{- \frac{\pi}{2}}^{\frac{\pi}{2}}{\frac{2 \cdot (2x)^{\frac{1}{2} + \frac{1}{\ln x}}}{\frac{1}{2}} \cdot \frac {d\varphi}{\ln x}} \ll
				\frac{\sqrt{x}}{\ln x}.$$
		Further,
		$$|F(s)| \le |\zeta(s)|^{\frac{1}{4}}(\ln x)^{\frac{19}{244}}|G(s)| \ll (\ln x)^{\frac{19}{244}}\left|\zeta(\sigma+it)\right|^{\frac{1}{4}}.$$
		Hence
		$$|j_2| = \left| \text{p.v.} \frac{1}{2 \pi i}\int_{\frac{1}{2} + iT}^{\frac{1}{2} + \frac{i}{\ln x}}{F(s) \cdot \frac {(x + h)^s - x^s}{s} ds} \right|
			\ll \int_{\frac{1}{\ln x}}^{T} {(\ln x)^{\frac{19}{244}} \cdot \left|\zeta\left(\tfrac{1}{2} + it\right)\right|^{\frac{1}{4}}} \sqrt{x} \frac{dt}{t+1} \ll $$
			$$\ll (\ln x)^{\frac{19}{244}} \sqrt{x} \int_{0}^{T}{\left|\zeta\left(\tfrac{1}{2}+it\right)\right|^{\frac{1}{4}} \frac{dt}{t+1}} =
			(\ln x)^{\frac{19}{244}} \sqrt{x} \sum_{\nu \ge 0}{\int_{T/2^\nu}^{T/2^{\nu +1}}{\frac{|\zeta(\frac{1}{2} +it)|^{\frac{1}{4}}}{t+1}dt}}.$$
		Denoting the summands in the last sum by $j(\nu)$ and taking $X=T\cdot2^{-\nu}$, by the Hölder inequality we get:
		$$j(\nu) \ll \frac{1}{X} \left( \int_{X}^{2X}{\left|\zeta\left(\tfrac{1}{2}+it\right)\right|^2dt}\right)^{\frac{1}{8}} X^{1-\frac{1}{8}} \ll
		\frac{1}{X} \left( X\ln X\right)^{\frac{1}{8}}X^{1-\frac{1}{8}} \ll
		(\ln X)^{\frac{1}{8}} \ll (\ln T)^{\frac{1}{8}}.$$
		Hence,
		$$\sum_{\nu \ge 0}{j(\nu)} \ll (\ln T)^{1+\frac{1}{8}} \ll (\ln T)^{\frac{9}{8}}.$$
		Then the upper bound for $j_2$ has the form
		$$|j_2| \ll (\ln x)^{\frac{19}{244}} \sqrt{x} (\ln x)^{\frac{9}{8}} = \sqrt{x} (\ln x)^{\frac{19}{244}+\frac{9}{8}} \ll \sqrt{x} (\ln x)^{\frac{27}{122}}.$$
		The integral $j_7$ is estimated as above.
		
		The main term arises from the calculation of $j_4$ and $j_5$. Let us define the entire function $w(s)$ by the relation
		$$\zeta(s) = \frac{w(s)}{s-1}$$
		and let $s = 1 - u + i\cdot0$, where $0 \le u \le \frac{1}{2}$. Then
		$$\sqrt[4]{\zeta(s)} = \frac {\sqrt[4]{w(s)}}{\sqrt[4]{-u+i \cdot 0}}.$$
		Since $-u+i\varepsilon \to u\cdot e^{\pi i}$ as $\varepsilon \to +0$, then
		$$\sqrt[4]{-u+i \cdot 0} = \sqrt[4]{u} e^{\frac{\pi i}{4}}, \;\; \sqrt[4]{\zeta(s)}=\frac{\sqrt[4]{w(\sigma)}}{\sqrt[4]{u}} e^{-\frac{\pi i}{4}}.$$
		Therefore, on the upper edge of the cut we have
		$$F(s) = \frac{\sqrt[4]{w(1-u)}}{(\zeta(2-2u))^{\frac{19}{244}}}\; G(1-u)\frac{e^{-\frac{\pi i}{4}}}{\sqrt[4]{u}} = \frac{\Pi(u)e^{-\frac{\pi i}{4}}}{\sqrt[4]{u}},$$
		where
		$$\Pi(u) = G(1-u) \frac{\sqrt[4]{w(1-u)}}{(\zeta(2-2u))^{\frac{19}{244}}}.$$
		
		Hence,
		$$j_4 = \frac{1}{2 \pi i} \int_{\frac{1}{2} + \frac{1}{\ln x} + i \cdot 0}^{1+i \cdot 0}{F(\sigma + i \cdot 0) \frac{(x+h)^s- x^s}{s} ds} = $$
		$$=\frac{1}{2\pi i} \int_{\frac{1}{2}+\frac{1}{\ln x}+i \cdot 0}^{1+i \cdot 0}{F(\sigma+i \cdot 0) \int_{0}^{h}{(x+u)^{s-1}du}\;ds} =$$
		$$= \frac{1}{2\pi i}\int_{x}^{x+h}{\int_{\frac{1}{2}+\frac{1}{\ln x}}^{1}{F(\sigma+i \cdot 0)y^{\sigma-1}d\sigma}\;dy}=
		\frac{e^{-\frac{\pi i}{4}}}{2\pi i} \int_{x}^{x+h}{{\int_{0}^{\frac{1}{2}-\frac{1}{\ln x}}{\frac{\Pi(u)y^{-u}}{\sqrt[4]{u}}}du}\;dy}.$$
		Suppose that $N \ge 0$ is fixed. Then
		$$\Pi(u)= \Pi_0+\Pi_1u + \Pi_2u^2 + \ldots + \Pi_N u^N + O_N(u^{N+1}),$$
		where
		$$\Pi_0=\Pi(0)=\frac{\sqrt[4]{w(1)}}{(\zeta(2))^{\frac{19}{244}}}.$$
		Thus, we have
		$$j_4 = \frac{e^{-\frac{\pi i}{4}}}{2\pi i} \int_{x}^{x+h}{\left(\sum_{0 \le n \le N}{\Pi_n \int_{0}^{\frac{1}{2}-\frac{1}{\ln x}}{\frac{u^{n}y^{-u}}{\sqrt[4]{u}}du}}+ O\left( J\right)\right)dy},$$
		where
		$$J = \int_{0}^{\frac{1}{2}-\frac{1}{\ln x}}{\frac{u^{N+1}y^{-u}}{\sqrt[4]{u}}du} \le  \frac{\Gamma\left(N+\frac{7}{4}\right)}{(\ln y)^{N+\frac{7}{4}}}.$$
		Using the estimate
		$$\int_{\lambda}^{\infty}{w^{k-\gamma}e^{-w}dw < ek!\;\lambda^{k-\gamma}e^{-\lambda}},$$
		where $\lambda>1$, $0 < \gamma <1$, $k \ge 1$, we easily get
		$$\int_{0}^{\frac{1}{2}-\frac{1}{\ln x}}{\frac{u^n y^{-u}}{\sqrt[4]{u}}du} =\frac{1}{(\ln y)^{n+\frac{3}{4}}} \int_{0}^{\ln y \left(\frac{1}{2}-\frac{1}{\ln x} \right)}{e^{-w} w^{n-\frac{1}{4}} dw}=$$
		$$=\frac{1}{(\ln y)^{n+\frac{3}{4}}} \left( \int_{0}^{\infty}{e^{-w} w^{n-\frac{1}{4}} dw} + \frac{\theta e n! (\ln y)^{n-\frac{1}{4}}}{\sqrt{y}}\right)=
		\frac{\Gamma\left( n+\frac{3}{4}\right)}{(\ln y)^{n+\frac{3}{4}}}+ \frac{\theta e n!}{\sqrt{y} \ln y}.$$
		Therefore,
		$$j_4 = \frac{e^{-\frac{\pi i}{4}}}{2\pi i} \int_{0}^{h}{\sum_{0 \le n \le N}{\frac{\Pi_n \Gamma(n+\frac{3}{4})}{(\ln (x+u))^{n+\frac{3}{4}}}}du}+O\left( \frac{h}{\left(\ln x\right)^{N+\frac{7}{4}}}\right).$$
		Let
		$$\varphi(x) = \frac{1}{\left( \ln x\right)^{n+\frac{7}{4}}}.$$
		Then the Lagrange mean-value theorem yields
		$$\varphi(x+u) = \varphi(x) + u \varphi^{'}(x+\theta_1 u) = \frac{1}{\left( \ln x\right)^{n+\frac{7}{4}}} + \frac{\theta_2 h \left( n+\frac{7}{4}\right)}{x \left( \ln x\right)^{n+\frac{7}{4}}}.$$
		Thus we get
		$$j_4 = \frac{h e^{-\frac{\pi i}{4}}}{2\pi i} \sum_{0 \le n \le N}{\frac{\Pi_n \Gamma\left( n+\frac{7}{4}\right)}{(\ln x)^{n+\frac{7}{4}}}}+ O\left( \frac{h}{(\ln x)^{N+\frac{7}{4}}}\right) + O\left( \frac{h^2}{x} \frac{1}{(\ln x)^{N+\frac{7}{4}}}\right),$$
		$$j_5 = -\frac{h e^{\frac{\pi i}{4}}}{2\pi i} \sum_{0 \le n \le N}{\frac{\Pi_n \Gamma\left( n+\frac{7}{4}\right)}{(\ln x)^{n+\frac{7}{4}}}}+ O\left( \frac{h}{(\ln x)^{N+\frac{7}{4}}}\right) + O\left( \frac{h^2}{x} \frac{1}{(\ln x)^{N+\frac{7}{4}}}\right).$$
		
		Finally,
		$$-(j_4+j_5)=
		-\frac{h}{(\ln x)^{\frac{3}{4}}} \left( \sum_{0 \le n \le N}{\frac{(-1)^n \Pi_n}{\Gamma \left( \frac{1}{4}-n\right) (\ln x)^n}} + O\left( \frac{1}{(\ln x)^{N+1}}\right) + O \left( \frac{h}{x(\ln x)^{N+1}}\right)\right).$$
		
		It remains to estimate the sum
		$$\sum_{|\gamma|<T}{j_{\rho}},\;j_{\rho}=I_1(\rho)+I_2(\rho),\;\text{where}\;\;\rho = \beta +i\gamma.$$
		
		\begin{figure}[tbh]
			\begin{center}
				\begin{tikzpicture}

					\draw [->](0,3.5) -- (1.25,3.5);
					\draw (1.25,3.5) -- (2.5,3.5);
					\draw (0,2) -- (1.25,2);
					\draw [<-](1.25,2) -- (2.5,2);

					\draw [->](0,2.75) -- (5,2.75);

					\draw (0,4.1) -- (0,3.5);
					\draw (0,1) -- (0,2);

					\fill[black] (3.4,2.75) circle(0.05);

					\draw[black] (2.5,3.5) arc (140:0:1.18);
					\draw[black] (2.5,2) arc (220:360:1.18);

					\draw (1.25,3.5) node [above]  {$I_1(\rho)$};
					\draw (1.25,2) node [below]  {$I_2(\rho)$};
					\draw (3.4,2.75) node [below]  {$\rho=\beta+i\gamma$};
					\draw (5,2.75) node [above]  {$\sigma$};
					\draw (0,1.2) node [right]  {$\scriptstyle{\Res s = \frac{1}{2}}$};

			\end{tikzpicture}
		\end{center}
		\end{figure}
		Since
		$$\left| \frac{(x+h)^s-x^s}{s}\right| = \left| \int_{0}^{h}{(x+u)^{s-1}du}\right| \ll  \int_{0}^{h}{(x+u)^{\sigma-1}du} \ll h x^{\sigma-1},$$
		then
		$$I_1(\rho) \ll \int_{\frac{1}{2}}^{\beta}{(\ln x)^{\frac{19}{244}} |\zeta(\sigma+i\gamma)|^{\frac{1}{4}} h x^{\sigma-1} d\sigma} \ll
		\frac{h}{x} (\ln x)^{\frac{19}{244}} \int_{\frac{1}{2}}^{\beta}{x^{\sigma}  |\zeta(\sigma+i\gamma)|^{\frac{1}{4}} d\sigma}$$
		and the same estimate is valid for $I_2(\rho)$.
		Hence,
		$$|j_\rho| \ll \int_{\frac{1}{2}}^{\beta}{h x^{\sigma-1} (\ln x)^{\frac{19}{244}} T^{\frac{c(1-\sigma)}{4}} (\ln x)^{\frac{1}{4}} d\sigma} \ll
		h (\ln x)^{\frac{19}{244}+\frac{1}{4}} \int_{\frac{1}{2}}^{\beta}{\left( \frac{T^{\frac{c}{4}}}{x}\right)^{1-\sigma} d\sigma} \ll$$
		$$\ll h (\ln x)^{\frac{19}{244}+\frac{1}{4}} \int_{\frac{1}{2}}^{1}{g(\rho,\sigma) \left( \frac{T^{\frac{c}{4}}}{x}\right)^{1-\sigma} d\sigma},$$
		where
		$$g(\rho,\sigma) = \left\{
			\begin{aligned}
				1,\;\text{if}\;\sigma \le \beta,\\ 
				0,\;\text{if}\;\sigma > \beta.\\
			\end{aligned}
			\right.$$
		Applying lemma \ref{density}, we get
		$$\sum_{|\gamma|<T}{j_\rho} \ll h (\ln x)^{\frac{20}{61}} \int_{\frac{1}{2}}^{1-\varrho (T)}{\left( \sum_{|\gamma| < T}{g(\rho; \gamma)} \right) \left( \frac{T^{\frac{c}{4}}}{x}\right)^{1-\sigma} d\sigma} \ll$$
		$$\ll h (\ln x)^{\frac{20}{61}} \int_{\frac{1}{2}}^{1-\varrho (T)} {N(\sigma; T) \left( \frac{T^{\frac{c}{4}}}{x}\right)^{1-\sigma} d\sigma}\ll$$
		$$\ll h (\ln x)^{44+\frac{20}{61}} \int_{\frac{1}{2}}^{1-\varrho (T)} {\left( \frac{T^{\frac{c}{4}}}{x}\right)^{1-\sigma} T^{\frac{12}{5}(1-\sigma)} d\sigma} \ll$$
		$$\ll  h (\ln x)^{45} \int_{\frac{1}{2}}^{1-\varrho (T)} {\left( \frac{T^{\frac{12}{5}+\frac{c}{4}}}{x}\right)^{1-\sigma} d\sigma} \ll
		 h (\ln x)^{45} \left( \frac{T^{\frac{12}{5} + \frac{c}{4}}}{x}\right)^{1-\varrho (T)}.$$
		Let $D(x)=e^{C_1(\ln x)^{0.8}}.$
		Choosing $T$ from the equation,
		$$T^{\frac{12}{5} + \frac{c}{4}}=x D^{-1}(x)$$
		we easily conclude that
		$$T=x^{\left(\frac{c}{4} + \frac{12}{5}\right)^{-1}} D(x)^{-\left(\frac{c}{4} + \frac{12}{5}\right) ^{-1}}.$$
		
		Obviously, the formula for $S_2$ is asymptotic if
		$$h \gg \frac{x}{T} (\ln x)^2 = x^{\alpha_2} e^{C_2 (\ln x)^{0.8}},$$	
		where
		$$\ = 1-\frac{205}{508}=\frac{303}{508}.$$
	\end{subsection}

	\begin{subsection}{The mean value of the function $2^{-\omega(n)}$ on the short interval}
		Suppose that $\sigma = \Res s > 1$ and let
		$$F(s) = \sum_{n=1}^{\infty}{\frac{1}{2^{\omega(n)}} \cdot n^{-s}}.$$
		This series converges absolutely, since 
		$$|F(s)| \le \sum_{n=1}^{\infty}{\left| \frac{1}{2^{\omega(n)}}\right| \cdot n^{-\sigma}} \le \frac{1}{2} \sum_{n=1}^{\infty}{n^{-\sigma}} = 
		\frac{1}{2} \left( 1 + \int_{1}^{\infty}{\frac{du}{u^\sigma}}\right) = \frac{1}{2} \left( 1 + \frac{1}{\sigma-1}\right).$$
		
		Setting $a_n = 2^{-\omega(n)}, A(n) \equiv 1, b = 1 +\displaystyle{\frac{1}{\ln x}}, \alpha =1$ in lemma \ref{perron}, we get
		$$S_3=S(x,h;f_3)=I+O(R),$$
		where
		$$I=\frac{1}{2\pi i}\int_{b-iT}^{b+iT}{F(s) \frac{(x+h)^s-x^s}{s}ds},\;\;\; R=\frac{x^b}{T(b-1)}+\frac{xA(2x)\ln x}{T} \ll \frac{x \ln x}{T}.$$
		Further, $F(s) = \prod_p{F_p(s)}$, where
		$$F_p(s)=1 + \frac{1}{2^{\omega(p)} p^{s}} + \frac{1}{2^{\omega(p^2)}p^{2s}} + \ldots = 1+\frac{1}{2 p^s} + \frac{1}{2 p^{2s}} + \ldots$$
		Writing $F_p(s)$ in the form
		$$F_p(s)=\left( 1 - \frac{1}{p^s}\right)^{-\frac{1}{2}} \left( 1 - \frac{1}{p^{2s}}\right)^{-\frac{1}{8}} G_p(s),$$
		we obtain
		$$F(s)=\frac{(\zeta(s))^{\frac{1}{2}}}{(\zeta(2s))^{\frac{1}{8}}} G(s),$$
		where
		$$G(s)=\prod_p{G_p(s)}=\prod_p{\left( 1 - \frac{1}{p^s}\right)^{\frac{1}{2}} \left( 1 - \frac{1}{p^{2s}}\right)^{\frac{1}{8}} \left(1+\frac{1}{2(p^s-1)}\right)}.$$
		Further
		$$G_p(s)=\left( 1 - \frac{1}{p^s}\right)^{-\frac{1}{2}} \left( 1 - \frac{1}{p^{2s}}\right)^{\frac{1}{8}} \left(1-\frac{1}{2p^s}\right).$$
		Now we continue the function $F(s)$ to the left of the line $\Res s =1$.
		Suppose that $\frac{1}{2} \le \sigma \le 1$. Then
		$$\ln G_p(s) = -\frac{1}{2} \ln \left( 1 - \frac{1}{p^s}\right)+
		\frac{1}{8} \ln \left( 1 - \frac{1}{p^{2s}}\right)+
		\ln \left( 1 - \frac{1}{2p^s}\right) = \sum_{n=1}^{\infty}{g_n p^{-ns}}.$$
		Thus, the coefficient for $p^{-ns}$ is equal to
		$$g_n = \begin{cases} 
		\frac{1}{n} \left( \frac{1}{4} - \frac{1}{2^n}\right), & \mbox{if } n\ge 2 \mbox{ is even,} \\
		\frac{1}{n} \left( \frac{1}{2} - \frac{1}{2^n}\right), & \mbox{if } n \ge 1 \mbox{ is odd.}
		\end{cases}
		$$
		Since $g_1 = 0$, $g_2 = 0$, then for all $n \ge 3$ we obtain $g_n \le \frac{1}{6}.$
		Next
		$$|\ln G_p(s)| \le \sum_{n=3}^{\infty}{g_n p^{-n\sigma}} \le 
		\frac{1}{6 p^{3\sigma}} \frac{1}{1- \frac{1}{p^\sigma}} \le \frac{1}{6p(\sqrt{p}-1)}.$$		
		
		Finally, for $\frac{1}{2} \le \sigma \le 1$ we get
		$$\left| \ln |G(s)|\right| \le \sum_{n=2}^{\infty}{\frac{1}{6p(\sqrt{p}-1)}} < C_3,$$
		$$-C_3 \le \ln |G(s)| \le C_3, \;\;\; e^{-C_3} \le |G(s)| \le e^{C_3}.$$
		
		Let $\Gamma$ be the boundary of the rectangle with the vertices $\frac{\Z 1}{\Z 2} \pm iT, b \pm iT$, where the zeros of $\zeta(s)$ of the form
		$\frac{\Z 1}{\Z 2}+i\gamma$, $|\gamma|<T$, are avoided by the semicircles of the infinitely small radius lying to the right of the line $\Res s = \frac{\Z 1}{\Z 2}$,
		the pole of $\zeta(2s)$ at the point $s = \frac{\Z 1}{\Z 2}$ is avoided by two arcs $\Gamma_1$ and $\Gamma_2$ with the radius $\frac{\Z 1}{\Z \ln x}$,
		and let a horizontal cut be drawn from the critical line inside this rectangle to each zero $\rho = \beta + i\gamma$, $\beta > \frac{\Z 1}{\Z 2}$, $|\gamma| < T$.
		Then the function $F(s)$ is analytic inside $\Gamma$.
		By the Cauchy residue theorem,
		
		$$j_0 = - \sum_{k=1}^{8}{j_k} - \sum_{\rho}{j_{\rho}} = - (j_4+j_5) - \sum_{k \neq 4,5}{j_k} - \sum_{\rho}{j_{\rho}}.$$
		
		\begin{figure}[tbh]
			\begin{center}
			\begin{picture}(250,267)
			\put(0,125){\vector(1,0){250}}
			\put(245,128){$\sigma$}
			
			\put(17,0){\vector(0,1){250}}
			\put(20,245){$t$}
			
			\put(17,125){\circle*{3}}
			\put(10,127){$0$}
			
			\put(17,10){\circle*{3}}
			\put(0,10){$-T$}
			
			\put(17,240){\circle*{3}}
			\put(7,235){$T$}
			
			\put(97,125){\circle*{3}}
			\put(100,115){$\frac{1}{2}$}
			
			\put(177,125){\circle*{3}}
			\put(165,115){$1$}
			\put(177,125){\oval(40,40)[r]}
			
			\put(97,145){\oval(20,20)[tr]}
			\put(97,105){\oval(20,20)[br]}
			
			\put(127,55){\circle*{3}}
			\put(127,55){\oval(20,20)[r]}
			
			\put(227,125){\circle*{3}}
			\put(230,115){$b$}
			
			\put(227,10){\vector(0,1){57}}
			\put(227,67){\line(0,1){58}}
			\put(230,67){$j_0$}
			
			\put(227,125){\vector(0,1){57}}
			\put(227,182){\line(0,1){58}}
			\put(230,182){$j_0$}
			
			\put(227,240){\vector(-1,0){65}}
			\put(162,240){\line(-1,0){65}}
			\put(162,247){$j_1$}
			
			\put(95,245){$\frac{1}{2}+iT$}
			\put(215,245){$b+iT$}
			
			\put(95,0){$\frac{1}{2}+iT$}
			\put(215,0){$b-iT$}
			
			\put(177,105){\vector(-1,0){40}}
			\put(137,105){\line(-1,0){30}}
			\put(137,95){$j_5$}
			
			\put(97,10){\vector(1,0){65}}
			\put(162,10){\line(1,0){65}}
			\put(162,0){$j_8$}
			
			\put(107,145){\vector(1,0){40}}
			\put(147,145){\line(1,0){30}}
			\put(137,152){$j_4$}
			
			\put(80,182){\vector(0,1){58}}
			\put(80,182){\vector(0,-1){57}}
			\put(68,182){$j_2$}
			
			\put(80,67){\vector(0,1){58}}
			\put(80,67){\vector(0,-1){57}}
			\put(68,67){$j_7$}
			
			\put(102,157){$j_3$}
			\put(102,88){$j_6$}
			
			\put(97,95){\line(0,-1){30}}
			\put(97,45){\vector(0,-1){18}}
			\put(97,27){\line(0,-1){17}}
			
			\put(97,240){\line(0,-1){10}}
			\put(97,220){\line(0,-1){20}}
			
			\put(97,190){\oval(20,20)[r]}
			\put(97,190){\circle*{3}}
			\put(97,180){\line(0,-1){25}}
			
			
			\put(97,65){\line(1,0){30}}
			\put(105,70){$I_1(\rho)$}
			
			\put(97,45){\line(1,0){30}}
			\put(105,35){$I_2(\rho)$}
			
			\put(127,55){\circle*{3}}
			\put(127,55){\oval(20,20)[r]}

			\put(97,230){\line(1,0){50}}
			\put(97,220){\line(1,0){50}}

			\put(147,225){\circle*{3}}
			\put(147,225){\oval(10,10)[r]}
			
			\multiput(17,240)(10,0){8}{\line(1,0){5}}
			\multiput(17,10)(10,0){8}{\line(1,0){5}}
			
			\multiput(60,150)(10,0){5}{\line(1,0){5}}
			\multiput(60,100)(10,0){5}{\line(1,0){5}}
			
			\put(48,150){$\Gamma_1$}
			\put(48,100){$\Gamma_2$}
			
			\put(100,150){\vector(1,0){5}}
			\put(100,100){\vector(1,0){5}}
			
		\end{picture}
			\end{center}
		\end{figure}
		
		By lemma \ref{T^()lnT},
		$$F(s) \ll T^{\frac{c(1-\sigma)}{2}} (\ln T)^{\frac{1}{8}}.$$
		Then
		$$|j_1|=\left| \frac{1}{2 \pi i} \int_{b + iT}^{\frac{1}{2} + iT}{F(s) \frac {(x + h)^s - x^s}{s} ds} \right| \ll
			\frac{1}{T} \int_{\frac{1}{2}}^{b}{T^{\frac{c(1 - \sigma)}{2}} \cdot (\ln x)^{\frac{1}{8}} x^{\sigma} d\sigma} \ll $$
		$$\ll \frac{x}{T} \int _{\frac{1}{2}}^{b}{\frac{x^{\sigma-1}}{T^{\frac{c(\sigma-1)}{2}}} (\ln x)^{\frac{1}{8}}} d\sigma \ll
		\frac{x}{T} \int _{\frac{1}{2}}^{b}{\left(\frac{x}{T^{\frac{c}{2}}} \right)^{\sigma-1} (\ln x)^{\frac{1}{8}}} d\sigma \ll
			\frac{x}{T} (\ln x)^{\frac{1}{8}}.$$
		The similar estimate is valid for $j_8$.
		
		\begin{figure}[tbh]
			\begin{center}
				\begin{picture}(125,150)
			\put(0,75){\vector(1,0){125}}
			\put(121,77){$\sigma$}
			
			\put(10,75){\circle*{3}}
			\put(9,64){$\frac{1}{2}$}
			
			\put(10,85){\oval(40,40)[tr]}
			\put(10,65){\oval(40,40)[br]}
			
			\put(30,85){\vector(1,0){30}}
			\put(60,85){\line(1,0){30}}
			
			\put(10,105){\line(0,1){10}}
			\put(10,35){\line(0,1){10}}
			
			\put(90,65){\vector(-1,0){30}}
			\put(30,65){\line(1,0){30}}
			
			\put(90,75){\oval(20,20)[r]}
			
			\put(25,105){$j_3$}
			\put(25,45){$j_6$}
			
			\put(0,95){$\Gamma_1$}
			\put(0,50){$\Gamma_2$}
			
			\put(12,97){\vector(1,0){13}}
			\put(12,52){\vector(1,0){13}}
			
		\end{picture}
			\end{center}
		\end{figure}
		
		By the estimations from theorems \ref{theorem1}, \ref{theorem2}, on $\Gamma_1, \Gamma_2$ we have:
		$$|\zeta(s)| \le 3.2, \left| \zeta(2s)\right| > 0.4 \ln x.$$
		Hence
		$$|F(s)|  \le \left| G(s) \right| \frac {{3.2}^{\frac{1}{2}}}{(0.4 \cdot \ln x)^{\frac{1}{8}}} < C.$$
		Therefore,
		$$|j_3 + j_6| \le \frac {1}{2\pi} \int_{\Gamma_1 \cup \Gamma_2}{|F(s)| \left| \frac {(x + h)^s - x^s}{s} \right| ds} \le
				\frac {C}{2\pi} \int_{- \frac{\pi}{2}}^{\frac{\pi}{2}}{\frac{2 \cdot (2x)^{\frac{1}{2} + \frac{1}{\ln x}}}{\frac{1}{2}} \cdot \frac {d\varphi}{\ln x}} \ll
				\frac{\sqrt{x}}{\ln x}.$$
		Further,
		$$|F(s)| \le |\zeta(s)|^{\frac{1}{2}}(\ln x)^{\frac{1}{8}}|G(s)| \ll (\ln x)^{\frac{1}{8}}\left|\zeta(\sigma+it)\right|^{\frac{1}{2}}.$$
		Hence
		$$|j_2| = \left| \text{p.v.} \frac{1}{2 \pi i}\int_{\frac{1}{2} + iT}^{\frac{1}{2} + \frac{i}{\ln x}}{F(s) \cdot \frac {(x + h)^s - x^s}{s} ds} \right|
			\ll \int_{\frac{1}{\ln x}}^{T} {(\ln x)^{\frac{1}{8}} \cdot \left|\zeta\left(\tfrac{1}{2} + it\right)\right|^{\frac{1}{2}}} \sqrt{x} \frac{dt}{t+1} \ll $$
			$$\ll (\ln x)^{\frac{1}{8}} \sqrt{x} \int_{0}^{T}{\left|\zeta\left(\tfrac{1}{2}+it\right)\right|^{\frac{1}{2}} \frac{dt}{t+1}} =
			(\ln x)^{\frac{1}{8}} \sqrt{x} \sum_{\nu \ge 0}{\int_{T/2^\nu}^{T/2^{\nu +1}}{\frac{|\zeta(\frac{1}{2} +it)|^{\frac{1}{2}}}{t+1}dt}}.$$
		Denoting the summands in the last sum by $j(\nu)$ and taking $X=T\cdot2^{-\nu}$, by the Hölder inequality we get:
		$$j(\nu) \ll \frac{1}{X} \left( \int_{X}^{2X}{\left|\zeta\left(\tfrac{1}{2}+it\right)\right|^2dt}\right)^{\frac{1}{4}} X^{1-\frac{1}{4}} \ll
		\frac{1}{X} \left( X\ln X\right)^{\frac{1}{4}}X^{1-\frac{1}{4}} \ll
		(\ln X)^{\frac{1}{4}} \ll (\ln T)^{\frac{1}{4}}.$$
		Hence,
		$$\sum_{\nu \ge 0}{j(\nu)} \ll (\ln T)^{1+\frac{1}{4}} \ll (\ln T)^{\frac{5}{4}}.$$
		Then the upper bound for $j_2$ has the form
		$$|j_2| \ll (\ln x)^{\frac{1}{8}} \sqrt{x} (\ln x)^{\frac{5}{4}} = \sqrt{x} (\ln x)^{\frac{1}{8}+\frac{5}{4}} \ll \sqrt{x} (\ln x)^{\frac{11}{8}}.$$
		The integral $j_7$ is estimated as above.
		
		The main term arises from the calculation of $j_4$ and $j_5$. Let us define the entire function $w(s)$ by the relation
		$$\zeta(s) = \frac{w(s)}{s-1}$$
		and let $s = 1 - u + i\cdot0$, where $0 \le u \le \frac{1}{2}$. Then
		$$\sqrt{\zeta(s)} = \frac {\sqrt{w(s)}}{\sqrt{-u+i \cdot 0}}.$$
		Since $-u+i\varepsilon \to u\cdot e^{\pi i}$ as $\varepsilon \to +0$, then
		$$\sqrt{-u+i \cdot 0} = \sqrt{u} e^{\frac{\pi i}{2}}, \;\; \sqrt{\zeta(s)}=\frac{\sqrt{w(\sigma)}}{\sqrt{u}} e^{-\frac{\pi i}{2}}.$$
		Therefore, on the upper edge of the cut we have
		$$F(s) = \frac{\sqrt{w(1-u)}}{(\zeta(2-2u))^{\frac{1}{8}}}\; G(1-u)\frac{e^{-\frac{\pi i}{2}}}{\sqrt{u}} = \frac{\Pi(u)e^{-\frac{\pi i}{2}}}{\sqrt{u}},$$
		where
		$$\Pi(u) = G(1-u) \frac{\sqrt{w(1-u)}}{(\zeta(2-2u))^{\frac{1}{8}}}.$$
		
		Hence,
		$$j_4 = \frac{1}{2 \pi i} \int_{\frac{1}{2} + \frac{1}{\ln x} + i \cdot 0}^{1+i \cdot 0}{F(\sigma + i \cdot 0) \frac{(x+h)^s- x^s}{s} ds} = $$
		$$=\frac{1}{2\pi i} \int_{\frac{1}{2}+\frac{1}{\ln x}+i \cdot 0}^{1+i \cdot 0}{F(\sigma+i \cdot 0) \int_{0}^{h}{(x+u)^{s-1}du}\;ds} =$$
		$$= \frac{1}{2\pi i}\int_{x}^{x+h}{\int_{\frac{1}{2}+\frac{1}{\ln x}}^{1}{F(\sigma+i \cdot 0)y^{\sigma-1}d\sigma}\;dy}=
		\frac{e^{-\frac{\pi i}{2}}}{2\pi i} \int_{x}^{x+h}{{\int_{0}^{\frac{1}{2}-\frac{1}{\ln x}}{\frac{\Pi(u)y^{-u}}{\sqrt{u}}}du}\;dy}.$$
		Suppose that $N \ge 0$ is fixed. Then
		$$\Pi(u)= \Pi_0+\Pi_1u + \Pi_2u^2 + \ldots + \Pi_N u^N + O_N(u^{N+1}),$$
		where
		$$\Pi_0=\Pi(0)=\frac{\sqrt{w(1)}}{(\zeta(2))^{\frac{1}{8}}}.$$
		Thus, we have
		$$j_4 = \frac{e^{-\frac{\pi i}{2}}}{2\pi i} \int_{x}^{x+h}{\left(\sum_{0 \le n \le N}{\Pi_n \int_{0}^{\frac{1}{2}-\frac{1}{\ln x}}{\frac{u^{n}y^{-u}}{\sqrt{u}}du}}+ O\left( J\right)\right)dy},$$
		where
		$$J = \int_{0}^{\frac{1}{2}-\frac{1}{\ln x}}{\frac{u^{N+1}y^{-u}}{\sqrt{u}}du} \le  \frac{\Gamma\left(N+\frac{3}{2}\right)}{(\ln y)^{N+\frac{3}{2}}}.$$
		Using the estimate
		$$\int_{\lambda}^{\infty}{w^{k-\gamma}e^{-w}dw < ek!\;\lambda^{k-\gamma}e^{-\lambda}},$$
		where $\lambda>1$, $0 < \gamma <1$, $k \ge 1$, we easily get
		$$\int_{0}^{\frac{1}{2}-\frac{1}{\ln x}}{\frac{u^n y^{-u}}{\sqrt{u}}du} =\frac{1}{(\ln y)^{n+\frac{1}{2}}} \int_{0}^{\ln y \left(\frac{1}{2}-\frac{1}{\ln x} \right)}{e^{-w} w^{n-\frac{1}{2}} dw}=$$
		$$=\frac{1}{(\ln y)^{n+\frac{1}{2}}} \left( \int_{0}^{\infty}{e^{-w} w^{n-\frac{1}{2}} dw} + \frac{\theta e n! (\ln y)^{n-\frac{1}{2}}}{\sqrt{y}}\right)=
		\frac{\Gamma\left( n+\frac{1}{2}\right)}{(\ln y)^{n+\frac{1}{2}}}+ \frac{\theta e n!}{\sqrt{y} \ln y}.$$
		Therefore,
		$$j_4 = \frac{e^{-\frac{\pi i}{2}}}{2\pi i} \int_{0}^{h}{\sum_{0 \le n \le N}{\frac{\Pi_n \Gamma(n+\frac{1}{2})}{(\ln (x+u))^{n+\frac{1}{2}}}}du}+O\left( \frac{h}{\left(\ln x\right)^{N+\frac{3}{2}}}\right).$$
		Let
		$$\varphi(x) = \frac{1}{\left( \ln x\right)^{n+\frac{3}{2}}}.$$
		Then the Lagrange mean-value theorem yields
		$$\varphi(x+u) = \varphi(x) + u \varphi^{'}(x+\theta_1 u) = \frac{1}{\left( \ln x\right)^{n+\frac{3}{2}}} + \frac{\theta_2 h \left( n+\frac{3}{2}\right)}{x \left( \ln x\right)^{n+\frac{3}{2}}}.$$
		Thus we get
		$$j_4 = \frac{h e^{-\frac{\pi i}{2}}}{2\pi i} \sum_{0 \le n \le N}{\frac{\Pi_n \Gamma\left( n+\frac{3}{2}\right)}{(\ln x)^{n+\frac{3}{2}}}}+ O\left( \frac{h}{(\ln x)^{N+\frac{3}{2}}}\right) + O\left( \frac{h^2}{x} \frac{1}{(\ln x)^{N+\frac{3}{2}}}\right),$$
		$$j_5 = -\frac{h e^{\frac{\pi i}{2}}}{2\pi i} \sum_{0 \le n \le N}{\frac{\Pi_n \Gamma\left( n+\frac{3}{2}\right)}{(\ln x)^{n+\frac{3}{2}}}}+ O\left( \frac{h}{(\ln x)^{N+\frac{3}{2}}}\right) + O\left( \frac{h^2}{x} \frac{1}{(\ln x)^{N+\frac{3}{2}}}\right).$$
		
		Finally,
		$$-(j_4+j_5)=
		-\frac{h}{(\ln x)^{\frac{1}{2}}} \left( \sum_{0 \le n \le N}{\frac{(-1)^n \Pi_n}{\Gamma \left( \frac{1}{2}-n\right) (\ln x)^n}} + O\left( \frac{1}{(\ln x)^{N+1}}\right) + O \left( \frac{h}{x(\ln x)^{N+1}}\right)\right).$$
		
		It remains to estimate the sum
		$$\sum_{|\gamma|<T}{j_{\rho}},\;j_{\rho}=I_1(\rho)+I_2(\rho),\;\text{where}\;\;\rho = \beta +i\gamma.
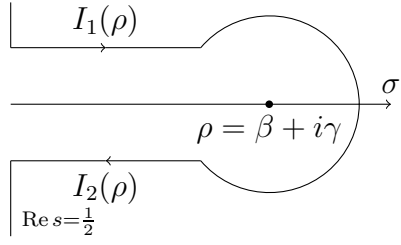
$$
		
		\begin{figure}[tbh]
			\begin{center}
				\begin{tikzpicture}

					\draw [->](0,3.5) -- (1.25,3.5);
					\draw (1.25,3.5) -- (2.5,3.5);
					\draw (0,2) -- (1.25,2);
					\draw [<-](1.25,2) -- (2.5,2);

					\draw [->](0,2.75) -- (5,2.75);

					\draw (0,4.1) -- (0,3.5);
					\draw (0,1) -- (0,2);

					\fill[black] (3.4,2.75) circle(0.05);

					\draw[black] (2.5,3.5) arc (140:0:1.18);
					\draw[black] (2.5,2) arc (220:360:1.18);

					\draw (1.25,3.5) node [above]  {$I_1(\rho)$};
					\draw (1.25,2) node [below]  {$I_2(\rho)$};
					\draw (3.4,2.75) node [below]  {$\rho=\beta+i\gamma$};
					\draw (5,2.75) node [above]  {$\sigma$};
					\draw (0,1.2) node [right]  {$\scriptstyle{\Res s = \frac{1}{2}}$};

			\end{tikzpicture}
		\end{center}
		\end{figure}
		Since
		$$\left| \frac{(x+h)^s-x^s}{s}\right| = \left| \int_{0}^{h}{(x+u)^{s-1}du}\right| \ll  \int_{0}^{h}{(x+u)^{\sigma-1}du} \ll h x^{\sigma-1},$$
		then
		$$I_1(\rho) \ll \int_{\frac{1}{2}}^{\beta}{(\ln x)^{\frac{1}{8}} |\zeta(\sigma+i\gamma)|^{\frac{1}{2}} h x^{\sigma-1} d\sigma} \ll
		\frac{h}{x} (\ln x)^{\frac{1}{8}} \int_{\frac{1}{2}}^{\beta}{x^{\sigma}  |\zeta(\sigma+i\gamma)|^{\frac{1}{2}} d\sigma}$$
		and the same estimate is valid for $I_2(\rho)$.
		Hence,
		$$|j_\rho| \ll \int_{\frac{1}{2}}^{\beta}{h x^{\sigma-1} (\ln x)^{\frac{1}{8}} T^{\frac{c(1-\sigma)}{2}} (\ln x)^{\frac{1}{2}} d\sigma} \ll
		h (\ln x)^{\frac{1}{8}+\frac{1}{2}} \int_{\frac{1}{2}}^{\beta}{\left( \frac{T^{\frac{c}{2}}}{x}\right)^{1-\sigma} d\sigma} \ll$$
		$$\ll h (\ln x)^{\frac{1}{8}+\frac{1}{2}} \int_{\frac{1}{2}}^{1}{g(\rho,\sigma) \left( \frac{T^{\frac{c}{2}}}{x}\right)^{1-\sigma} d\sigma},$$
		where
		$$g(\rho,\sigma) = \left\{
			\begin{aligned}
				1,\;\text{if}\;\sigma \le \beta,\\ 
				0,\;\text{if}\;\sigma > \beta.\\
			\end{aligned}
			\right.$$
		Applying lemma \ref{density}, we get
		$$\sum_{|\gamma|<T}{j_\rho} \ll h (\ln x)^{\frac{5}{8}} \int_{\frac{1}{2}}^{1-\varrho (T)}{\left( \sum_{|\gamma| < T}{g(\rho; \gamma)} \right) \left( \frac{T^{\frac{c}{2}}}{x}\right)^{1-\sigma} d\sigma} \ll$$
		$$\ll h (\ln x)^{\frac{5}{8}} \int_{\frac{1}{2}}^{1-\varrho (T)} {N(\sigma; T) \left( \frac{T^{\frac{c}{2}}}{x}\right)^{1-\sigma} d\sigma}\ll$$
		$$\ll h (\ln x)^{44+\frac{5}{8}} \int_{\frac{1}{2}}^{1-\varrho (T)} {\left( \frac{T^{\frac{c}{2}}}{x}\right)^{1-\sigma} T^{\frac{12}{5}(1-\sigma)} d\sigma} \ll$$
		$$\ll  h (\ln x)^{45} \int_{\frac{1}{2}}^{1-\varrho (T)} {\left( \frac{T^{\frac{12}{5}+\frac{c}{2}}}{x}\right)^{1-\sigma} d\sigma} \ll
		 h (\ln x)^{45} \left( \frac{T^{\frac{12}{5} + \frac{c}{2}}}{x}\right)^{1-\varrho (T)}.$$
		Let $D(x)=e^{C_1(\ln x)^{0.8}}.$
		Choosing $T$ from the equation,
		$$T^{\frac{12}{5} + \frac{c}{2}}=x D^{-1}(x)$$
		we easily conclude that
		$$T=x^{\left(\frac{c}{2} + \frac{12}{5}\right)^{-1}} D(x)^{-\left(\frac{c}{2} + \frac{12}{5}\right) ^{-1}}.$$
		
		Obviously, the formula for $S_3$ is asymptotic if
		$$h \gg \frac{x}{T} (\ln x)^2 = x^{\alpha_3} e^{C_2 (\ln x)^{0.8}},$$	
		where
		$$\alpha_3 = 1-\frac{205}{524}=\frac{319}{524}.$$
\end{subsection}
	
\begin{subsection}{The mean value of the function $2^{-\Omega(n)}$ on the short interval}
		Suppose that $\sigma = \Res s > 1$ and let
		$$F(s) = \sum_{n=1}^{\infty}{\frac{1}{2^{\Omega(n)}} \cdot n^{-s}}.$$
		This series converges absolutely, since 
		$$|F(s)| \le \sum_{n=1}^{\infty}{\left| \frac{1}{2^{\Omega(n)}}\right| \cdot n^{-\sigma}} \le \frac{1}{2} \sum_{n=1}^{\infty}{n^{-\sigma}} = 
		\frac{1}{2} \left( 1 + \int_{1}^{\infty}{\frac{du}{u^\sigma}}\right) = \frac{1}{2} \left( 1 + \frac{1}{\sigma-1}\right).$$
		
		Setting $a_n = 2^{-\Omega(n)}, A(n) \equiv 1, b = 1 +\displaystyle{\frac{1}{\ln x}}, \alpha =1$ in lemma \ref{perron}, we get
		$$S_4=S(x,h;f_4)=I+O(R),$$
		where
		$$I=\frac{1}{2\pi i}\int_{b-iT}^{b+iT}{F(s) \frac{(x+h)^s-x^s}{s}ds},\;\;\; R=\frac{x^b}{T(b-1)}+\frac{xA(2x)\ln x}{T} \ll \frac{x \ln x}{T}.$$
		Further, $F(s) = \prod_p{F_p(s)}$, where
		$$F_p(s)=1 + \frac{1}{2^{\Omega(p)} p^{s}} + \frac{1}{2^{\Omega(p^2)}p^{2s}} + \ldots = 1+\frac{1}{2 p^s} + \frac{1}{2^2 p^{2s}} + \ldots$$
		Writing $F_p(s)$ in the form
		$$F_p(s)=\left( 1 - \frac{1}{p^s}\right)^{-\frac{1}{2}} \left( 1 - \frac{1}{p^{2s}}\right)^{\frac{1}{8}} G_p(s),$$
		we obtain
		$$F(s)=\frac{(\zeta(s))^{\frac{1}{2}}}{(\zeta(2s))^{\frac{1}{8}}} G(s),$$
		where
		$$G(s)=\prod_p{G_p(s)}=\prod_p{\left( 1 - \frac{1}{p^s}\right)^{\frac{1}{2}} \left( 1 - \frac{1}{p^{2s}}\right)^{-\frac{1}{8}} \left(1-\frac{1}{2p^s}\right)^{-1}}.$$
		Now we continue the function $F(s)$ to the left of the line $\Res s =1$.
		Suppose that $\frac{1}{2} \le \sigma \le 1$. Then
		$$\ln G_p(s) = \frac{1}{2} \ln \left( 1 - \frac{1}{p^s}\right)-
		\frac{1}{8} \ln \left( 1 - \frac{1}{p^{2s}}\right)-
		\ln \left( 1 - \frac{1}{2p^s}\right) = \sum_{n=1}^{\infty}{g_n p^{-ns}}.$$
		Thus, the coefficient for $p^{-ns}$ is equal to
		$$g_n = \begin{cases} 
		\frac{1}{n} \left(\frac{1}{2^n}-\frac{1}{4}\right), & \mbox{if } n\ge 2 \mbox{ is even,} \\
		\frac{1}{n} \left( \frac{1}{2^n} - \frac{1}{2}\right), & \mbox{if } n \ge 1 \mbox{ is odd.}
		\end{cases}
		$$
		Since $g_1 = 0$, $g_2 = 0$, then for all $n \ge 3$ we obtain $g_n \le \frac{1}{6}.$
		Next
		$$|\ln G_p(s)| \le \sum_{n=3}^{\infty}{g_n p^{-n\sigma}} \le 
		\frac{1}{6 p^{3\sigma}} \frac{1}{1- \frac{1}{p^\sigma}} \le \frac{1}{6p(\sqrt{p}-1)}.$$		
		
		Finally, for $\frac{1}{2} \le \sigma \le 1$ we get
		$$\left| \ln |G(s)|\right| \le \sum_{n=2}^{\infty}{\frac{1}{6p(\sqrt{p}-1)}} < C_3,$$
		$$-C_3 \le \ln |G(s)| \le C_3, \;\;\; e^{-C_3} \le |G(s)| \le e^{C_3}.$$
		
		Let $\Gamma$ be the boundary of the rectangle with the vertices $\frac{\Z 1}{\Z 2} \pm iT, b \pm iT$, where the zeros of $\zeta(s)$ of the form
		$\frac{\Z 1}{\Z 2}+i\gamma$, $|\gamma|<T$, are avoided by the semicircles of the infinitely small radius lying to the right of the line $\Res s = \frac{\Z 1}{\Z 2}$,
		the pole of $\zeta(2s)$ at the point $s = \frac{\Z 1}{\Z 2}$ is avoided by two arcs $\Gamma_1$ and $\Gamma_2$ with the radius $\frac{\Z 1}{\Z \ln x}$,
		and let a horizontal cut be drawn from the critical line inside this rectangle to each zero $\rho = \beta + i\gamma$, $\beta > \frac{\Z 1}{\Z 2}$, $|\gamma| < T$.
		Then the function $F(s)$ is analytic inside $\Gamma$.
		By the Cauchy residue theorem,
		
		$$j_0 = - \sum_{k=1}^{8}{j_k} - \sum_{\rho}{j_{\rho}} = - (j_4+j_5) - \sum_{k \neq 4,5}{j_k} - \sum_{\rho}{j_{\rho}}.$$
		
		\begin{figure}[tbh]
			\begin{center}
			\begin{picture}(250,267)
			\put(0,125){\vector(1,0){250}}
			\put(245,128){$\sigma$}
			
			\put(17,0){\vector(0,1){250}}
			\put(20,245){$t$}
			
			\put(17,125){\circle*{3}}
			\put(10,127){$0$}
			
			\put(17,10){\circle*{3}}
			\put(0,10){$-T$}
			
			\put(17,240){\circle*{3}}
			\put(7,235){$T$}
			
			\put(97,125){\circle*{3}}
			\put(100,115){$\frac{1}{2}$}
			
			\put(177,125){\circle*{3}}
			\put(165,115){$1$}
			\put(177,125){\oval(40,40)[r]}
			
			\put(97,145){\oval(20,20)[tr]}
			\put(97,105){\oval(20,20)[br]}
			
			\put(127,55){\circle*{3}}
			\put(127,55){\oval(20,20)[r]}
			
			\put(227,125){\circle*{3}}
			\put(230,115){$b$}
			
			\put(227,10){\vector(0,1){57}}
			\put(227,67){\line(0,1){58}}
			\put(230,67){$j_0$}
			
			\put(227,125){\vector(0,1){57}}
			\put(227,182){\line(0,1){58}}
			\put(230,182){$j_0$}
			
			\put(227,240){\vector(-1,0){65}}
			\put(162,240){\line(-1,0){65}}
			\put(162,247){$j_1$}
			
			\put(95,245){$\frac{1}{2}+iT$}
			\put(215,245){$b+iT$}
			
			\put(95,0){$\frac{1}{2}+iT$}
			\put(215,0){$b-iT$}
			
			\put(177,105){\vector(-1,0){40}}
			\put(137,105){\line(-1,0){30}}
			\put(137,95){$j_5$}
			
			\put(97,10){\vector(1,0){65}}
			\put(162,10){\line(1,0){65}}
			\put(162,0){$j_8$}
			
			\put(107,145){\vector(1,0){40}}
			\put(147,145){\line(1,0){30}}
			\put(137,152){$j_4$}
			
			\put(80,182){\vector(0,1){58}}
			\put(80,182){\vector(0,-1){57}}
			\put(68,182){$j_2$}
			
			\put(80,67){\vector(0,1){58}}
			\put(80,67){\vector(0,-1){57}}
			\put(68,67){$j_7$}
			
			\put(102,157){$j_3$}
			\put(102,88){$j_6$}
			
			\put(97,95){\line(0,-1){30}}
			\put(97,45){\vector(0,-1){18}}
			\put(97,27){\line(0,-1){17}}
			
			\put(97,240){\line(0,-1){10}}
			\put(97,220){\line(0,-1){20}}
			
			\put(97,190){\oval(20,20)[r]}
			\put(97,190){\circle*{3}}
			\put(97,180){\line(0,-1){25}}
			
			
			\put(97,65){\line(1,0){30}}
			\put(105,70){$I_1(\rho)$}
			
			\put(97,45){\line(1,0){30}}
			\put(105,35){$I_2(\rho)$}
			
			\put(127,55){\circle*{3}}
			\put(127,55){\oval(20,20)[r]}

			\put(97,230){\line(1,0){50}}
			\put(97,220){\line(1,0){50}}

			\put(147,225){\circle*{3}}
			\put(147,225){\oval(10,10)[r]}
			
			\multiput(17,240)(10,0){8}{\line(1,0){5}}
			\multiput(17,10)(10,0){8}{\line(1,0){5}}
			
			\multiput(60,150)(10,0){5}{\line(1,0){5}}
			\multiput(60,100)(10,0){5}{\line(1,0){5}}
			
			\put(48,150){$\Gamma_1$}
			\put(48,100){$\Gamma_2$}
			
			\put(100,150){\vector(1,0){5}}
			\put(100,100){\vector(1,0){5}}
			
		\end{picture}
			\end{center}
		\end{figure}
		
		Since
		$$F(s)= \frac{(\zeta(s))^{\frac{1}{2}}}{(\zeta(2s))^{\frac{1}{8}}} G(s)$$
		is exactly the same as in theorem \ref{theorem3}, then all the previous estimates hold. Thus, we have
	
		$$|j_1| \ll \frac{x}{T} (\ln x)^{\frac{1}{8}},\quad |j_8| \ll \frac{x}{T} (\ln x)^{\frac{1}{8}},$$
		$$|j_3 + j_6| \ll \frac{\sqrt{x}}{\ln x},$$
		$$|j_2| \ll \sqrt{x} (\ln x)^{\frac{11}{8}}, \quad|j_7| \ll \sqrt{x} (\ln x)^{\frac{11}{8}},$$		
		$$-(j_4+j_5)=-\frac{h}{(\ln x)^{\frac{1}{2}}} \left( \sum_{0 \le n \le N}{\frac{(-1)^n \Pi_n}{\Gamma \left( \frac{1}{2}-n\right) (\ln x)^n}} + O\left( \frac{1}{(\ln x)^{N+1}}\right) + O \left( \frac{h}{x(\ln x)^{N+1}}\right)\right),
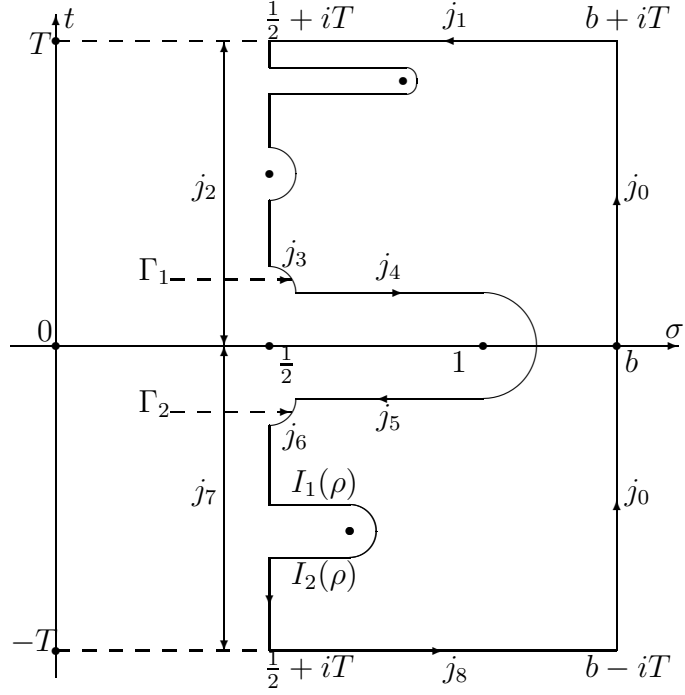
$$
		$$\sum_{|\gamma|<T}{j_\rho} \ll 
		 h (\ln x)^{45} \left( \frac{T^{\frac{12}{5} + \frac{c}{2}}}{x}\right)^{1-\varrho (T)}.$$
		Let $D(x)=e^{C_1(\ln x)^{0.8}}.$
		Choosing $T$ from the equation,
		$$T^{\frac{12}{5} + \frac{c}{2}}=x D^{-1}(x)$$
		we easily conclude that
		$$T=x^{\left(\frac{c}{2} + \frac{12}{5}\right)^{-1}} D(x)^{-\left(\frac{c}{2} + \frac{12}{5}\right) ^{-1}}.$$
		
		Obviously, the formula for $S_3$ is asymptotic if
		$$h \gg \frac{x}{T} (\ln x)^2 = x^{\alpha_4} e^{C_2 (\ln x)^{0.8}},$$	
		where
		$$\alpha_4 = 1-\frac{205}{524}=\frac{319}{524}.$$
\end{subsection}	

\end{section}


\begin{thebibliography}{}
\bibitem{me}{A. Sedunova, On the mean values of some multiplicative functions on the short interval, {\em Preprint: arxiv.org, arXiv:1302.0471, 2013.}}
	\bibitem{rmnj}{S. Ramanujan, Some Formulae in the Analythic Theory of Numbers. {\em Mess. Math.\/}, 45(1916), pp. 81-84.}
	\bibitem{wilson}{B.M. Wilson, Proofs of some formulae enunciated by Ramanujan. {\em Proc. London Math. Soc.\/}, 2(21), 1922, pp. 235-255.}
	\bibitem{mukan}{A.T. Mukanova, An asymptotic formula for the mean value of the V. I. Arnold function. {In Russian.} {\em Vestnik Moskovskogo Universiteta\/}, 2008, Vol. 63, No. 2, pp. 51–53.}
	\bibitem{erdos}{P.T. Bateman, P. Erdös, C. Pomerance, E.G.Straus, The arithmetic mean of the divisors of an integer. {\em Analytic number theory, Lecture Notes in Math.\/},1981, pp. 197-220.}
	\bibitem{vinogradov}{I.M. Vinogradov, Elements of Number Theory, Courier Dover Publications, 2003.}
	\bibitem{vorkara}{A.A. Karatsuba, S.M. Voronin, The Riemann Zeta-Function. Walter de Gruyter, 1992.}
	\bibitem{huxley}{M.N. Huxley, On the difference between consecutive primes. {\em Invent. math.\/}, 15, 1972, pp. 164-170.}
	\bibitem{rmch}{K. Ramachandra, On the number of Goldbach numbers in small intervals. {\em J. Indian Math. Soc.\/}, 37, 1973, pp. 157-170.}
	\bibitem{lavrik}{A.F. Lavrik, The approximate functional equation for Dirichlet $L$-functions. {In Russian.} {\em Tr. Mosk. Mat. Obs.\/}, 18, MSU, M., 1968, pp. 91–104.}
	\bibitem{titchmarch}{E.C. Titchmarsh, The Theory of the Riemann Zeta-Function. Oxford University Press Inc., 1951.}
	\bibitem{kar}{A.A. Karatsuba, Basic Analytic Number Theory. Springer-Verlag, 1992.}
	\bibitem{fcttheory}{E.C. Titchmarsh, The Theory of Functions. Oxford University Press, 1939.}
	\bibitem{rockafellar}{R. Rockafellar, Convex analysis. Princeton University Press, 1970.}
	\bibitem{huxley2}{M.N.Huxley,Exponential sums and the Riemann Zeta Function V. {\em Proc. London Math. Soc\/}, 90, 2005, 1-41.}
\end{thebibliography}
\end{document}